\documentclass[twocolumn]{autart}   
\usepackage{amsmath}
\usepackage{color}
\usepackage{amsfonts,amssymb}
\usepackage{subfigure}
\usepackage{graphicx}   
\usepackage{algorithm}
\usepackage{algorithmicx}  
\usepackage{algpseudocode} 
\usepackage{makecell}
\begin{document}

\begin{frontmatter}

\title{Nonlinear Set Membership Filter with State Estimation Constraints via Consensus-ADMM\thanksref{footnoteinfo}} 

\thanks[footnoteinfo]{This work was supported in part by the NSFC under Grant 61673282 and Grant U1836103.}

\author[scu]{Xiaowei Li}\ead{lxwstu@outlook.com},
\author[scu]{Xuqi Zhang}\ead{zxqcc@stu.scu.edu.cn}, 
\author[scu]{Zhiguo Wang}\ead{wangzhiguo@scu.edu.cn},  
\author[scu]{Xiaojing Shen}\ead{shenxj@scu.edu.cn}  

\address[scu]{School of Mathematics, Sichuan University, Chengdu 610064, China}                
\begin{abstract}
This paper considers the state estimation problem for nonlinear dynamic systems with unknown but bounded noises. Set membership filter (SMF) is a popular algorithm to solve this problem. In the set membership setting, we investigate the filter problem where the state estimation requires to be constrained by a linear or nonlinear equality. We propose a consensus alternating direction method of multipliers (ADMM) based SMF algorithm for nonlinear dynamic systems. To deal with the difficulty of nonlinearity, instead of linearizing the nonlinear system, a semi-infinite programming (SIP) approach is used to transform the nonlinear system into a linear one, which allows us to obtain a more accurate estimation ellipsoid. For the solution of the SIP, an ADMM algorithm is proposed to handle the state estimation constraints, and each iteration of the algorithm can be solved efficiently. Finally, the proposed filter is applied to typical numerical examples to demonstrate its effectiveness.
\end{abstract}

\end{frontmatter}

\section{Introduction}
State estimation problems of dynamic systems occur in many applications, such as robot localization \cite{saeedi2016multiple}, target tracking \cite{7950937}, machine learning \cite{6841003}, fault diagnosis \cite{simon2006kalman}, navigation \cite{Alouani1993}, etc. Since it is difficult to measure the system state directly, it is crucial to estimate the state from noisy sensor measurements, imperfect system models, and physical constraints. In the setting of stochastic noises, Bayesian filtering has been extensively researched. Specifically, Kalman filter \cite{kalman1960new} is a special case of Bayesian filtering under the linear, quadratic, and Gaussian conditions \cite{chen2003bayesian}, and it is the minimum-variance recursive state estimator. For the nonlinear dynamic systems, many modified Kalman-based filters \cite{simon2010,ristic2003beyond,Zhu2002Multisensor} are proposed to approximate the optimal state estimation. When the statistical properties of the process and the measurement noises can be obtained precisely, these filters work well and are extensively applied to target tracking, navigation, and other applications.

Unlike the Bayesian strategy, the set membership method aims to provide guaranteed enclosures for the system state in the presence of unknown but bounded uncertainty effects that do not require the assumption of knowledge of its stochastic properties \cite{polyak2004ellipsoidal}. SMF is first considered by Schweppe \cite{schweppe1968recursive} in the late 1960s. Since then, the idea of SMF has been widely studied \cite{el2001robust,shen2011minimizing,wang2019ellipsoidal,9435052,durieu2001multi,wang2020set}. It is worth noting that most studies on the SMF have focused on the state estimation problems for unconstrained linear dynamic systems. Since the estimation given by SMF is generally a regular state bounding region while a nonlinear system maps a regular set to an irregular one. In addition, the shape of the state bounding set is further affected if there are additional constraints on the dynamic systems. Thus, extending SMF to nonlinear dynamic systems is difficult, especially when online usage needs to be guaranteed. Some studies on nonlinear dynamic systems are presented in \cite{wang2019ellipsoidal,9435052,chen2018nonlinear,rego2021set}. For example, in the extended set membership filter (ESMF) \cite{scholte2003nonlinear} and nonlinear set membership filter (NSMF) \cite{calafiore2005reliable}, the nonlinear functions in the dynamic system are linearized around the current state estimate, after which the higher-order remainder terms are bounded in different ways. Thus, the nonlinear dynamic system is transformed into a linear dynamic system by treating the bounds of the higher-order remainder terms as noises. The dual set membership filter (DSMF) proposed in \cite{9435052} is also a nonlinear set membership filter. Unlike ESMF and NSMF, DSMF directly finds the bounding ellipsoid of the nonlinearly transformed state bounding ellipsoid via dual theory and SIP, rather than linearizing the nonlinear dynamic system. Therefore, a tighter state bounding ellipsoid can be obtained by DSMF, while it only considers the unconstrained case. 

In many practical scenarios, the state estimation is usually expected to be localized in some particular subspace \cite{xiong2003blind,zhou2006linear,duan2016modeling}. For example, in a blind multiuser detection problem, it is necessary to ensure that the target value associated with the desired user is located in the specific searching space \cite{xiong2003blind}. In navigation and localization, state estimation of vehicles in three dimensional space is expected to be localized on the ground, while the localization of an operating aircraft is in the air \cite{kirubarajan2000ground}. Therefore, it is significant to consider state estimators with various estimation constraints. Previous studies on set membership filters mainly focus on state constraints, such as \cite{yang2009set1,yang2009set,li2020set} and the references therein. However, even if the state constraint information is used in the filtering process to ensure that the state bounding ellipsoids take into account the state constraint, the state estimates that are given by the ellipsoidal SMF (the center of the state bounding ellipsoids) generally do not satisfy the constraint. Motivated by this, the paper considers the nonlinear SMF problem with estimation constraints.

In this paper, we investigate the state estimation problem for nonlinear dynamic systems with unknown but bounded noises and state estimation constraints. The main contributions of the paper are as follows:
\begin{itemize}
  \item We propose a consensus-ADMM-based SMF algorithm to estimate the states of nonlinear dynamic systems with a general form equality estimation constraint. The centers of the state bounding ellipsoids are guaranteed to satisfy the estimation constraint, which makes up for the deficiency of SMF with the state constraint.
  \item To deal with the difficulty of nonlinearity, we use the SIP approach to transform the nonlinear dynamic system into a linear one rather than linearizing the nonlinear dynamic system, which allows us to determine a more accurate estimation ellipsoid.
  \item We provided a weak convergence result of the proposed consensus-ADMM algorithm. When we solve the discretized SIP problems, each iteration of the consensus-ADMM can be solved efficiently. Especially, for the linear estimation constraint, we can obtain an analytical solution for the state update.
\end{itemize}
Two numerical examples in the simulations show the effectiveness of the proposed filter. 

The rest of this paper is organized as follows. In Section \ref{Preliminaries}, the preliminaries and problem formulation are given. A nonlinear SMF with state estimation constraints for computing the state bounding ellipsoid is developed in Section \ref{iii}. In Section \ref{sectionadmm}, A consensus-ADMM-based method for solving the SIP problems is given. Simulations and conclusions are given in Section \ref{v} and \ref{vi}, respectively.

\textbf{Notations:} $\oplus$ represents the Minkowski sum. $G_{m,n}$ is the projection matrix from dimension $n$ to $m$. $A^\dagger$ represents the Moore-Penrose generalized inverse of matrix $A$. The notation $A \succeq 0$ means $A$ is positive semidefinite and $A \succ 0$ means $A$ is positive definite. The notation $\mathcal N(A)$ denotes the nullspace of $A$. $I_{\theta\times\theta}$ denotes the identity matrix with dimension $\theta$. In this paper, $J(E)$ is either $tr(E)$ or $\log \det(E)$.

{\setlength\abovedisplayskip{5pt}\setlength\belowdisplayskip{5pt}\section{Preliminaries and problem formulation}\label{Preliminaries}
\subsection{Preliminaries}

The basic idea of the SMF is to obtain a feasible solution set containing the true state according to the dynamic system, measurements, bounded noises, and other information about the state \cite{polyak2004ellipsoidal,el2001robust}. Specifically, suppose the initial state $x_0$ belongs to a given set ${\mathcal Y}_{0}$. The objective of the SMF is to determine a minimum volume state bounding set ${\mathcal Y}_{k+1}$ based on ${\mathcal Y}_{k}$, the dynamic system, measurements, bounded noises, and other information at time step $k+1$. 

The SMF is subdivided into the prediction and the measurement update step. In the prediction step at time $k+1$, a predicted state bounding set ${\mathcal Y}_{k+1|k}$ is determined as the minimum volume set according to the state bounding set ${\mathcal Y}_{k}$, the state transform equation, the process noise information, and other information. In the measurement update step at time $k+1$, a state bounding set ${\mathcal Y}_{k+1}$ is determined as the minimum volume set according to the predicted state bounding set ${\mathcal Y}_{k+1|k}$, the measurement equation, the measurements, the measurement noise information, and other information.

In general, two basic set representations are used in the SMF: polytopes of various types (e.g., general polytopes \cite{avis1997good}, boxes \cite{li2020generalized}, zonotopes \cite{girard2005reachability}, parallelotopes \cite{kostousova2001control}, and rectangular polytopes \cite{stursberg2003efficient}) and ellipsoid \cite{kurzhanskiy2007ellipsoidal}. Polytopes can give arbitrarily close approximations to any convex set, but the number of vertices can grow prohibitively large \cite{avis1997good}. The complexity of the ellipsoidal representation is quadratic in the dimension of the set since its size, shape, and location can be determined uniquely by its center and shape matrix. In addition, it is possible to single out a set approximating ellipsoid that is optimal to some given criterion or a combination of them \cite{kurzhanskiy2007ellipsoidal}. Therefore, we focus on the ellipsoid approach in this paper, and some preliminaries of the ellipsoid are provided as follows.

\textbf{Definition 1.} A set ${\mathcal E}\subset \mathbb{R}^n$ satisfying the following form is called an ellipsoid:
\begin{align*}
{\mathcal E}(e,E)=\{x\in \mathbb{R}^n: (x-e)^TE^{-1}(x-e)\leq1\},
\end{align*}
where $e\in \mathbb{R}^n$ and $E\succ 0$ are the center and the shape matrix of $\mathcal E$, respectively. 

The ``size” of an ellipsoid ${\mathcal E}(e,E)$ can be expressed as a function $J(E)$ of the shape matrix $E$ \cite{el2001robust}. In this paper, we choose the most common ways to measure the ``size" of an ellipsoid, which are the trace function ($J(E) = tr(E)$, which corresponds to the sum of squares of semiaxes lengths of $\mathcal E(e,E)$) and logdet function ($J(E) =\log \det(E)$, which corresponds to the volume of $\mathcal E(e,E)$), respectively \cite{durieu2001multi}. The problem of finding the smallest ellipsoidal outer approximation of a bounded set $\mathcal X\subset\mathbb{R}^n$ on the trace or logdet criterion can be described as:
\begin{align*}
&\mathop{\min\limits_{e,E}J(E)},~s.t.~\mathcal E(e,E) \supset{\mathcal X}.
\end{align*}\textbf{Definition 2.} The Minkowski sum of two ellipsoids ${\mathcal E}_1$ and ${\mathcal E}_2$ is defined as:
\begin{align*}
{\mathcal E}_1\oplus{\mathcal E}_2=\{x_1+x_2:x_1\in{\mathcal E}_1,x_2\in{\mathcal E}_2\}.
\end{align*}
\textbf{Definition 3.} The intersection of two ellipsoids ${\mathcal E}_1$ and ${\mathcal E}_2$ is defined as:
\begin{align*}
{\mathcal E}_1\cap{\mathcal E}_2=\{x:x\in{\mathcal E}_1,x\in{\mathcal E}_2\}.
\end{align*}
Next, we briefly introduce the consensus-ADMM method used in this paper \cite{boyd2011distributed,huang2016consensus}. Consider the problem:
\begin{align}
\min_{x}\sum^{\omega}_{i=1}\Pi_i(x),\label{prelimadmm}
\end{align}
where $x\in \mathbb{R}^n$, the objective function is the sum of the objective terms $\Pi_i(x):\mathbb{R}^n\to\mathbb{R}\cup\{+\infty\},~i=1,...,\omega$. Problem \eqref{prelimadmm} can be reformulated as a consensus form by introducing $\omega$ auxiliary variables $z_i\in \mathbb{R}^n,~i=1,...,\omega$:
\begin{align*}
&\min_{x,z_i}\sum^{\omega}_{i=1}\Pi_i(z_i),~s.t.~x=z_i,~\forall i=1,...,\omega,
\end{align*}
then the consensus-ADMM iterates at iteration $t$ for this problem are:
\begin{align*}
z_i^{t+1}:= &\mathop{\arg\min}\limits_{z_i}
\Pi_i(z_i)-(\lambda_i^{t})^Tz_i+(\rho/2)\|x^{t}-z_i\|^2,\nonumber\\
x^{t+1}:= &\mathop{\arg\min}\limits_{x}\sum^{\omega}_{i=1}((\lambda_i^{t})^Tx+(\rho/2)\|x-z_i^{t+1}\|^2),\nonumber\\
\lambda_i^{t+1}:= &\lambda_i^{t}+\rho(x^{t+1}-z^{t+1}_i), i=1,...,\omega,
\end{align*}where $\rho > 0$ is the penalty parameter, $\lambda_i$ is the Lagrange multiplier, and $x^{t}$ and $\lambda_i^{t}$ are the updates obtained in the $t$-th iteration and $z_i^{t+1}$, $x^{t+1}$, and $\lambda_i^{t+1}$ are the updates obtained in the $t+1$-th iteration, respectively. In the consensus-ADMM iterates, $z_i$, $x$, and $\lambda_i$ are updated in an alternating fashion, and the $z_i$-minimization step is independently for each $i$. More details on the consensus-ADMM method can be seen in \cite{boyd2011distributed}. The consensus-ADMM algorithm is a very intuitive algorithm with several advantages \cite{boyd2011distributed,huang2016consensus}. First, the update of the variables $z_i$ can be implemented in parallel since they are independent of each other. Next, each update can be efficient for the high-dimensional but sparse datasets by splitting the cost functions. Moreover, the consensus-ADMM method converges if the problem is convex.

\subsection{Problem formulation}\label{ii}
Consider the nonlinear dynamic system:
\begin{align}
x_{k+1}&=f_k(x_k)+w_k,\label{d1}\\
y_k&=h_k(x_k)+v_k,\label{d1yk}
\end{align} 
where $k$ denotes the time step, $x_k\in \mathbb{R}^n $ is the system state, \(y_k\in \mathbb{R}^m\) is the measurement, \(f_k(x_k)\) and \(h_k(x_k)\) are nonlinear uniformly continuous and differentiable process function and measurement function of $x_k$, respectively. The process noise \(w_k\in \mathbb{R}^n\) and the measurement noise \(v_k\in \mathbb{R}^m\) are assumed to be restricted in given ellipsoidal sets $W_k={\mathcal E}(0,Q_k)$ and $V_k={\mathcal E}(0,R_k)$, respectively. In many practical scenarios, the state estimation $\hat{x}_k$ is expected to be localized in specific subspace \cite{zhou2006linear}. In particular, we consider the following equality state estimation constraint:\begin{align}
g_k(\hat{x}_k)= 0,\label{ineq_1}
\end{align}
where $g_k(\hat{x}_k)$ is a continuously differentiable function.

Suppose the initial state $x_0$ belongs to ${\mathcal E}_{0}(\hat{x}_{0},P_{0})$. The objective of this paper is to design a SMF that determines a minimum volume state bounding ellipsoid ${\mathcal E}_{k+1}(\hat{x}_{k+1},P_{k+1})$ based on ${\mathcal E}_{k}(\hat{x}_{k},P_{k})$, the dynamic system, measurements, bounded noises, and constraint at time step $k+1$.

\section{Nonlinear set membership filter with state estimation constraints}\label{iii}

This section presents a nonlinear SMF for computing the state bounding ellipsoid with the center satisfying the constraint (\ref{ineq_1}) at each time step, where the prediction and the update ellipsoid can be obtained by solving SIP problems without linearizing the nonlinear system.

\subsection{Prediction step}\label{sectiona}

In this subsection, we consider the prediction step of the proposed nonlinear SMF. Let ${\mathcal E}_{k}(\hat{x}_{k},P_{k})$ and ${\mathcal E}_{k+1|k}(\hat{x}_{k+1|k},P_{k+1|k})$ denote the state bounding ellipsoid of $x_k$ and the predicted ellipsoid of $x_{k+1}$, respectively. According to the uniform continuity of $f_k$, we have that ${\mathcal F}_k=\{f_k(x_k): x_k\in{\mathcal E}_{k}\}$ is a bounded compact set. Therefore, a minimum volume ellipsoid ${\mathcal E}_{f_k}(\hat{x}_{f_k},P_{f_k})$ containing ${\mathcal F}_k$ can be obtained by solving: 
\begin{align}
&\min_{\hat{x}_{f_k},P_{f_k}} \quad J(P_{f_k})\nonumber \\
&s.t.\quad (x-\hat{x}_{f_k})^TP_{f_k}^{-1}(x-\hat{x}_{f_k})\leq1,~\forall x\in{\mathcal F}_k.\label{ADMM_0}
\end{align}
From \eqref{d1}, we have $x_{k+1}\in{\mathcal E}_{f_k}\oplus W_k$. Thus, the predicted state bounding ellipsoid at time $k+1$ can be obtained by solving:
\begin{align}
&\mathop{\min}\limits_{\hat x_{k+1|k},P_{k+1|k}}J(P_{k+1|k})\nonumber\\
&s.t.\quad \mathcal E_{k+1|k}(\hat x_{k+1|k},P_{k+1|k}) \supset{\mathcal E}_{f_k}\oplus W_k.
 \label{pre1}
\end{align}
\textbf{Remark 1.} Theorem 4.2 in \cite{durieu2001multi} and (\ref{pre1}) show that the shape matrix of ${\mathcal E}_{k+1|k}$ has the following form:
\begin{align}
P_{k+1|k}(\tau_k)=(1+\tau_k^{-1})P_{f_k}+(1+\tau_k)Q_k,\label{p1}
\end{align} 
where $\tau_k>0$. The selection of $\tau_k$ determines the property of ${\mathcal E}_{k+1|k}$, and the problem (\ref{pre1}) can be reformulated as $min_{\tau_k>0}J(P_{k+1|k}(\tau_k))$. Especially, if $J(P_{k+1|k})=tr(P_{k+1|k})$, then the optimal solution to this optimization problem is $\tau_k=\sqrt{\frac{tr(P_{f_k})}{tr(Q_k)}}$.

The following lemma shows the condition under which the state prediction satisfies the estimation constraint.

\textbf{Lemma 1.} The state prediction $\hat{x}_{k+1|k}$ satisfies the constraint (\ref{ineq_1}) if and only if $\hat{x}_{f_k}$ satisfies (\ref{ineq_1}).

\textbf{Proof:} See Appendix \ref{app0}. \hfill{$\square $}

Based on the optimization problem (\ref{ADMM_0}) and Lemma 1, in order to make the predicted state located in the constraint subspace, the problem of determining the predicted state bounding ellipsoid has been transformed into:
\begin{align}
&\min_{\hat{x}_{f_k},P_{f_k}} \quad J(P_{f_k})\nonumber \\
&s.t.\quad (x-\hat{x}_{f_k})^TP_{f_k}^{-1}(x-\hat{x}_{f_k})\leq1,~\forall x\in{\mathcal F}_k,\nonumber \\
&\qquad \qquad g_{k+1}(\hat{x}_{f_k})= 0.\label{ADMM_1}
\end{align}
The optimization problem (\ref{ADMM_1}) is a SIP problem that includes two variables and infinite inequality constraints, and it is not jointly convex in the two variables. Therefore, the problem (\ref{ADMM_1}) is difficult to solve, and we provide a consensus-ADMM-based method to handle it in Section \ref{sectionadmm}.

\subsection{Measurement update step}\label{sectionb}
Subsequently, we consider the measurement update step of the proposed filter. Based on the predicted ellipsoid and the measurement at time $k+1$, we seek the minimum volume state bounding ellipsoid ${\mathcal E}_{k+1}(\hat{x}_{k+1},P_{k+1})$ with the center $\hat{x}_{k+1}$ satisfying the estimation constraint (\ref{ineq_1}).

For the convenience of the analysis, we assume that there exists a continuous inverse function $h_{k+1}^{-1}$ for the nonlinear function $h_{k+1}$. Thus, the measurement function is reformulated as:
\begin{align}
G_{m,n}x_{k+1}=h_{k+1}^{-1}(y_{k+1}-v_{k+1}).\label{pro1}
\end{align} 

\textbf{Remark 2.} Equation \eqref{pro1} uses the information about the inverse function of the measurement function to estimate the state instead of linearizing the nonlinear system, thus leading to more accurate state estimates \cite{9435052}. In general, the assumption that there exists a continuous inverse function for $h_{k+1}$ can be very restrictive. Here, it can be satisfied in some practical applications \cite{9435052}. For example, the nonlinear measurement functions in two or three dimensional radar systems \cite{pu2018optimal}, which are the foundation for the measurement systems of many sensors \cite{julier2004unscented}. In addition, this assumption can be relaxed in some special and important fields \cite{9435052}. For example, the simultaneous localization and mapping (SLAM) problem in mobile robot localization \cite{thrun2002probabilistic}.

Since the measurement noise $v_{k+1}\in V_{k+1}$ and the function $h_{k+1}^{-1}$ is continuous, the set ${\mathcal S}_{k+1}=\{h_{k+1}^{-1}(y_{k+1}-v_{k+1}):v_{k+1}\in{V}_{k+1}\}$ on the right-hand side of (\ref{pro1}) is a compact set. Thus, an ellipsoid ${\mathcal E}_{h_{k+1}}(\hat{x}_{h_{k+1}},P_{h_{k+1}})$ containing $G_{m,n}x_{k+1}$ can be obtained by solving:
\begin{align}
&\min_{\hat{x}_{h_{k+1}},P_{h_{k+1}}}~J(P_{h_{k+1}})\nonumber \\
&s.t.~(x-\hat{x}_{h_{k+1}})^TP_{h_{k+1}}^{-1}(x-\hat{x}_{h_{k+1}})\leq1,\nonumber \\
&\qquad\qquad\qquad \forall x\in{\mathcal S}_{k+1},\label{ADMM_21}
\end{align}which is also a SIP problem and can be seen as a simplified form of (\ref{ADMM_1}).

The update step is to find a state bounding ellipsoid ${\mathcal E}_{k+1}$ such that ${\mathcal E}_{h_{k+1}}\cap {\mathcal E}_{k+1|k}\subseteq{\mathcal E}_{k+1}$, where ${\mathcal E}_{h_{k+1}}$ is the measurement ellipsoid and ${\mathcal E}_{k+1|k}$ is the predicted ellipsoid. By denoting ${\mathcal D}_{k+1}=\{x:x\in{\mathcal E}_{h_{k+1}}\cap {\mathcal E}_{k+1|k}\}$, ${\mathcal E}_{k+1}$ can be derived by solving:
\begin{align}
&\min_{\hat{x}_{{k+1}},P_{{k+1}}} \quad J(P_{{k+1}})\nonumber \\
&s.t. \quad (x-\hat{x}_{{k+1}})^TP_{{k+1}}^{-1}(x-\hat{x}_{{k+1}})\leq1,~\forall x\in{\mathcal D}_{k+1},\nonumber \\
&\qquad\qquad g_{k+1}(\hat{x}_{k+1})=0.\label{ADMM_31}
\end{align}
So far, the nonlinear set membership state estimation problem with estimation constraints has been transformed to solve the SIP problems (\ref{ADMM_1}), (\ref{ADMM_21}), and (\ref{ADMM_31}). Notably, the nonlinear functions in the dynamic system do not need to be linearized. In summary, the proposed nonlinear SMF with state estimation constraints is given in Algorithm \ref{al3}.

\textbf{Remark 3.} The proposed method in Algorithm \ref{al3} is subdivided into two phases: 1) A prediction step determined by the state equation \eqref{d1}, process noise, and constraint \eqref{ineq_1}. 2) A measurement update step determined by the measurement equation \eqref{d1yk}, measurements, measurement noise, and constraint \eqref{ineq_1}. Authors in \cite{ghaoui1999worst} studied the pure state prediction problem (without measurement information). By using the constraint information in the prediction step, the prediction steps can be concatenated if some measurements are missing \cite{polyak2004ellipsoidal}. Thus, the proposed method can also be used for the state prediction problem by solving the problem \eqref{ADMM_1} to obtain the predicted ellipsoid. When the predicted ellipsoid is not required in practice, we can determine the predicted ellipsoid by solving the problem \eqref{ADMM_0}. In Section \ref{v}, we compare the performance of Algorithm \ref{al3} between the prediction step with and without using the constraint information.

\section{Solving the SIP problems}\label{sectionadmm}
In this section, we focus on how to solve the SIP problems (\ref{ADMM_1}), (\ref{ADMM_21}), and (\ref{ADMM_31}). These problems can be unified as:
\begin{align}
&\min_{\hat{x},P}~ \log \det(P)\nonumber \\
&s.t.~ (r-\hat{x})^TP^{-1}(r-\hat{x})\leq 1,~\forall r\in{\mathcal T},~g(\hat{x})=0,\label{ADMM1}
\end{align}
where $\hat x\in \mathbb{R}^n$, $P\succ 0$, $\mathcal T\subset\mathbb{R}^n$ is a compact set, and $g(\hat{x})$ is a continuously differentiable function. If the objective function of (\ref{ADMM1}) is the trace function, similar results and algorithms can be obtained.
\begin{algorithm}[htbp]
\caption{Set membership filter with state estimation constraints}\label{al3}
\begin{algorithmic}
\Require The nonlinear functions $f_k$ and $h_{k+1}$, the measurement $y_{k+1}$, the ellipsoidal sets $W_k$ and $V_{k+1}$, the estimation constraint function $g_{k+1}$, the initial ellipsoid ${\mathcal E}_{0}(\hat{x}_{0},P_{0})$;
\Ensure
\For {each k=1:T}
\State Solve the optimization problem (\ref{ADMM_1}) (see Algorithm \ref{al4});
\State Calculate the shape matrix of the predicted ellipsoid by (\ref{p1});
\State Solve the optimization problem (\ref{ADMM_21}) to obtain $\hat{x}_{h_{k+1}}$ and $P_{h_{k+1}}$ (see Algorithm \ref{al4});
\State Calculate the measurement update ellipsoid by solving (\ref{ADMM_31}) (see Algorithm \ref{al4});
\EndFor\\
\Return {$\hat{x}_{k+1}$ and $P_{k+1}$};
\end{algorithmic}
\end{algorithm}

In general, the SIP is NP-hard \cite{ben2002tractable,wu2017sdr}. Since it is challenging to solve the SIP, we consider approximating it by some relaxation methods. One approach for solving the SIP problem is to minimize its objective function subject to only a finite subset of the infinite set of constraints \cite{hettich1986implementation,reemtsen1998semi,polak1992rate,lopez2007semi}, namely the discretization method. The discretization method has several characteristics \cite{liu2016sample}. First, the discretization method is distributionally robust since it works for general SIP. Next, by replacing the infinite constraints with finitely many simple constraints, the sample approximation technique significantly simplifies SIP. Last but not least, solving the sample approximation problem returns a solution to the original SIP with guaranteed performance. The discretized optimization problem of (\ref{ADMM1}) can be expressed as:\begin{align}
&\min_{\hat{x},P} \quad \log \det(P)\nonumber \\
&s.t. \quad (r_i-\hat{x})^TP^{-1}(r_i-\hat{x})\leq 1,~g(\hat{x})=0,~i\in I,\label{ADMM2}
\end{align}
where $r_i\in{\mathcal T}$, $I=\{1,2,...,s\}$, and $s$ is the sampling number.

\textbf{Remark 4.} The discretization method is a sample approximation scheme that randomly samples a sufficient number of constraints from the infinite constraints. In \cite{calafiore2005uncertain}, the authors provide an efficient and explicit bound on the measure of the original constraints that are possibly violated by the randomized solution, where the volume rapidly decreases to zero as the sample number increases. The authors in \cite{9435052} show that it only requires to sample from the boundary of the ellipsoid set for the nonlinear measurement functions in two or three dimensional radar systems. Thus, it can reduce a lot of redundant samples so that the computational time can be significantly decreased. Specifically, for the problems \eqref{ADMM_1} and \eqref{ADMM_21}, we choose $f_k(\bar r_i^1)$ and $h_{k+1}^{-1}(y_{k+1}-\bar r_i^2)$ as the samples $r_i$, respectively, where $\bar r_i^1$ and $\bar r_i^2$ are sampled from the ellipsoids ${\mathcal E_k}$ and $V_{k+1}$, respectively. For the problem \eqref{ADMM_31}, the set ${\mathcal D}_{k+1}$ is the intersection of two ellipsoids and the samples $r_i$ can be obtained by using the Accept-Reject method or Monte Carlo methods \cite{rabiei2021intersection}.

To solve problem (\ref{ADMM2}), we introduce the auxiliary variables $z_i,~i=1,...,s$, and then problem (\ref{ADMM2}) can be transformed into a consensus form \cite{boyd2011distributed}:
\begin{align}
&\min_{\hat{x},P,z_i} \quad \log \det(P)\nonumber \\
&s.t. \quad (r_i-z_i)^TP^{-1}(r_i-z_i)\leq 1,\nonumber \\
&\qquad\quad \hat{x}=z_i,~g(\hat{x})=0,~i\in I.\label{ADMM_41}
\end{align}
The corresponding augmented Lagrangian of (\ref{ADMM_41}) is given by:
\begin{align}
&L_\rho(P,z_i,\hat{x},\lambda_i)=\log \det(P)+\nonumber\\&\quad\qquad\sum_{i\in I}\lambda_i^T(\hat{x}-z_i)+\frac{\rho}{2}\sum_{i\in I}\|\hat{x}-z_i\|^2.
\label{}
\end{align}
The dual problem corresponding to the primal problem is as follows:
\begin{align}
&\max_{\lambda_i}\quad \min_{\tiny{\begin{matrix}_{\hat{x},P,z_i}\\(r_i-z_i)^TP^{-1}(r_i-z_i)\leq 1,\\ g(\hat{x})=0, i\in I\end{matrix}}} \qquad L_\rho(P,z_i,\hat{x},\lambda_i).\label{dual11}
\end{align}
By applying the consensus-ADMM method, we obtain the following iterations:
\begin{align}
z_i,P\leftarrow &\mathop{\arg\min}\limits_{z_i,P}
L_\rho(P,z_i,\hat{x},\lambda_i),\nonumber\\
&\quad s.t. \quad(r_i-z_i)^TP^{-1}(r_i-z_i)\leq 1, i\in I,\label{pz1}\\
\hat{x}\leftarrow &\mathop{\arg\min}\limits_{\hat{x}}
\sum_{i\in I}\lambda_i^T(\hat{x}-z_i)+\frac{\rho}{2}\sum_{i\in I}\|\hat{x}-z_i\|^2,\nonumber\\
&\quad s.t. \quad g(\hat{x})=0,\label{x0}\end{align}
\begin{align}
\lambda_i\leftarrow &\lambda_i+\rho(\hat{x}-z_i), i\in I\label{ADMM11111}.
\end{align}
It is well known that the ADMM algorithm converges under mild conditions when solving convex problems \cite{huang2016consensus,eckstein1992douglas}. Unfortunately, in nonconvex problems, convergence results cannot always be guaranteed. Next, we provide a weak convergence result for the proposed ADMM iterates of the optimization problem (\ref{ADMM2}).

\textbf{Theorem 1.} Denote $P^{t}$, $z_i^{t}$ and $\hat x^{t}$ as the updates obtained in the $t$-th
iteration of (\ref{pz1}) - (\ref{ADMM11111}). Assume that $z_i^{t}$ are well-defined for all $t$ and $i$, and that:
\begin{align}
&lim_{t\rightarrow\infty}(z_i^{t}-\hat x^{t})=0,\quad \forall i\in I,\label{prop1}\\
&lim_{t\rightarrow\infty}(\hat x^{t+1}-\hat x^{t})=0,\label{prop2}
\end{align}
then any limit point of $\{P^{t},\hat x^{t}\}$ is a KKT point of (\ref{ADMM2}).

\textbf{Proof:} See Appendix \ref{appa}. \hfill{$\square $}

The term ``well-defined” means existence and being uniquely defined, and this is a common assumption in convergence analysis that is rarely violated in practice \cite{huang2016consensus,xu2012alternating}.

\subsection{Update of $P$ and $z_i$}
In this subsection, we focus on the update of $P$ and $z_i$. Note that the objective function and constraints of the problem (\ref{pz1}) are convex in $P^{-1}$ and $z_i$, respectively. Nevertheless, (\ref{pz1}) is not convex in $P^{-1}$ and $z_i$ jointly because of the cross terms in the constraints. Thus, to simplify the iterations of $P$ and $z_i$, the block coordinate descent technique can be used \cite{bertsekas1997nonlinear}. Specifically, the first block of the ADMM iterates can take the following form:
\begin{align}
P\leftarrow &\mathop{\arg\min}\limits_{P}
\log \det(P)\nonumber\\
&\quad s.t.\quad (r_i-z_i)^TP^{-1}(r_i-z_i)\leq 1,~i\in I, \label{P111}\\
z_i\leftarrow &\mathop{\arg\min}\limits_{z_i}\lambda_i^T(\hat{x}-z_i)+\frac{\rho}{2}\|\hat{x}-z_i\|^2_2,\nonumber\\
&\quad s.t.\quad (r_i-z_i)^TP^{-1}(r_i-z_i)\leq 1.\label{ADMM_pz}
\end{align}
Note that both optimization problems (\ref{P111}) and (\ref{ADMM_pz}) are convex and the update of each variable $z_i$ are parallelizable. The following theorem gives a first-order optimal condition of (\ref{P111}).

\textbf{Theorem 2.} The optimal solution of (\ref{P111}) is
\begin{align}
{P}^*=n\sum_{i\in I}\mu_i^*(r_i-z_i)(r_i-z_i)^T,\label{pb}
\end{align}
where $\mu^{*}=(\mu_1^*,\mu_2^*,...,\mu_s^*)$ is the optimal solution of the dual problem:
\begin{align}
\max_\mu\quad &\log \det(\sum_{i\in I}\mu_i(r_i-z_i)(r_i-z_i)^T)\nonumber\\
s.t. \quad&\sum_{i\in I}\mu_i=1, \mu\geq 0.\label{P12}
\end{align} 
\textbf{Proof:} See Appendix \ref{appb}. \hfill{$\square $}

\textbf{Remark 5.} A projection-free first-order method for solving the problem (\ref{P12}) is the Frank-Wolfe (FW) algorithm \cite{todd2016minimum}, which iterates mainly by maximizing the linear Taylor approximation around the current solution $\mu$ on the unit simplex. Let $\varXi(\mu)=\log \det(\sum_{i\in I}\mu_i(r_i-z_i)(r_i-z_i)^T)$, the details of the FW algorithm are shown in Algorithm \ref{al1}. The global convergence analysis of the FW method has been given in the literature, e.g., \cite{todd2016minimum,damla2008linear,ahipacsaouglu2013modified}, and the literature shows that the number of iterations to obtain an $\epsilon$-approximately optimal solution is at most $O(1/\epsilon)$. Therefore, the maximum number of iterations $K$ can be determined according to the given accuracy $\epsilon$.

\begin{algorithm}[htbp]
\caption{FW algorithm for (\ref{P12})}
\begin{algorithmic}
\Require Choose $\mu^0$ satisfies $\sum_{i\in I}\mu^0_i=1$, the maximum number of iterations $K$;
\Ensure
\For {each $t=1:K$}
\State Compute $\varrho^{t}= \mathop{\arg\min}\limits_{\sum_{i\in I}\varrho_i=1} \varrho^T\partial \varXi(\mu^{t})$;
\State Obtain the optimal update step $\kappa^{t}$;
\State update $\mu^{t+1}=\mu^{t}+\kappa^{t}(\varrho^{t}-\mu^{t})$;
\EndFor\\
\Return {$\mu^{t+1}$};
\end{algorithmic}\label{al1}
\end{algorithm}

Next, we consider the update of $z_i$. The optimization problem (\ref{ADMM_pz}) is a quadratic constraint quadratic programming (QCQP-1) problem, and its solution satisfies the following theorem.

\textbf{Theorem 3.} $z_i^*$ is the optimal solution of optimization problem (\ref{ADMM_pz}) if and only if:
\begin{align}
z_i^*=(I+\eta_i^*P^{-1})^{-1}(\hat{x}-r_i+\frac{1}{\rho}\lambda_i)+r_i,
\end{align}
where $\eta_i^*=max\{0, \phi_i^*\}$, $\phi_i^*$ is the largest solution of $g(\phi_i)=1$, and $g(\phi_i)$ is defined as:
\begin{align}
g(\phi_i)=(\hat{x}-r_i+&\frac{1}{\rho}\lambda_i)^T(I-\phi_iP^{-1})^{-1}P^{-1}\nonumber\\&\cdot(I-\phi_iP^{-1})^{-1}(\hat{x}-r_i+\frac{1}{\rho}\lambda_i).\nonumber
\end{align}
\textbf{Proof:} See Appendix \ref{appc}. \hfill{$\square $}

\textbf{Remark 6.} The iteration of $z_i$ requires solving an equation $g(\phi_i)=1$. According to the proof of Theorem 3, $g(\phi_i)$ is a monotonically decreasing function which can be efficiently solved by many methods, e.g., the bisection method and Newton's method \cite{boyd2004convex}.

\subsection{Update of $\hat{x}$}
Now we consider the update of $\hat{x}$. For the optimization problem (\ref{ADMM_21}), the updated $\hat{x}$ can be obtained by solving:
\begin{align}
\min_{\hat{x}}\sum_{i\in I}\lambda_i^T(\hat{x}-z_i)+\frac{\rho}{2}\sum_{i\in I}\|\hat{x}-z_i\|^2,\label{x11}
\end{align}
where the optimal solution is $\hat{x}^*=\frac{1}{s}\sum_{i\in I}(z_i-\frac{1}{\rho}\lambda_i)$.

For the optimization problems (\ref{ADMM_1}) and (\ref{ADMM_31}), the updated $\hat{x}$ can be obtained by solving:
\begin{align}
\min_{\hat{x}}&\quad\sum_{i\in I}\lambda_i^T(\hat{x}-z_i)+\frac{\rho}{2}\sum_{i\in I}\|\hat{x}-z_i\|^2,\nonumber\\
 s.t.&\quad g(\hat{x})=0.\label{123}
\end{align}
The necessary optimality conditions of (\ref{123}) is \cite{boyd2004convex}:
\begin{align}
&\sum_{i\in I}\lambda_i^T+\rho\sum_{i\in I}(\hat{x}-z_i)+\nabla g(\hat{x})u=0,~g(\hat{x})=0,\label{1234}
\end{align}
where $u$ is the Lagrange multiplier vector. The nonlinear equations (\ref{1234}) are called the Lagrangian system and can be solved by first-order and second-order methods, such as Lagrangian and Newton’s methods \cite{bertsekas1997nonlinear}. Different algorithms to slove the problem (\ref{123}) can be selected according to the different form of the constraint function. Specifically, if $g$ is a linear function, then (\ref{123}) is a quadratic programming (QP) problem, and we have the following theorem.

\textbf{Theorem 4.} For a linear estimation constraint $g(\hat x)=C\hat x-c$, the update of $\hat{x}^*$ can be computed analytically and is given by:
\begin{align}
\hat{x}^*=C^\dagger c+UU^\dagger(d-C^\dagger c),\label{xxx}
\end{align} where $d=\frac{1}{s}\sum_{i\in I}(z_{i}-\frac{1}{\rho}\lambda_{i})$ and $U=I-C^\dagger C$ is an orthogonal projection matrix in $\mathcal N(C)$.

\textbf{Proof:} See Appendix \ref{appd}. \hfill{$\square $}

To sum up, the detailed consensus-ADMM algorithm for solving the problem (\ref{ADMM1}) is given in Algorithm \ref{al4}.
\begin{algorithm}
\caption{Consensus-ADMM algorithm for (\ref{ADMM1})}\label{al4}
\begin{algorithmic}
\Require The set $\mathcal T$, the number of samples $s$, penalty parameter $\rho>0$, tolerance $\epsilon>0$; Initialize $\hat x$, $z_i$ and $\lambda_i$; 
\Ensure 
\State Generate samples $r_1, r_2,...,r_s$ from set $\mathcal T$;
\Repeat 
\State Solve optimization problem (\ref{P12}) to get the dual variable $\mu$ (see Algorithm \ref{al1});
\State Calculate $P$ by equation (\ref{pb});
\For {each $i=1:s$}
\State Solve the QCQP-1 (\ref{ADMM_pz}) to get $z_i$;
\EndFor
\If{$g$ is a linear function}
\State Calculate $\hat{x}$ by (\ref{xxx});
\Else
\State Solve the Lagrangian system (\ref{1234}) to get $\hat{x}$;
\EndIf
\State Calculate $\lambda_i=\lambda_i+\rho(\hat{x}-z_i)$;
\Until the successive difference of $\hat x$ is smaller than a tolerance $\epsilon$, or the number of iterations reaches a given maximum;\\
\Return {$\hat{x}$ and $P$};
\end{algorithmic}
\end{algorithm}

\textbf{Remark 7.} When the equality constraints cannot hold exactly due to uncertainties, the perturbed equality constraints can be regarded as inequality constraints. By replacing the equality constraint in the SIP problem with the inequality constraint and solving an inequality constrained problem when updating $\hat x$, the proposed method can be extended to the inequality constraint case. For example, if the inequality constraint function is linear, then the proposed SMF algorithm needs to solve the linear inequality constrained SIP problems. Moreover, the problem of updating $\hat x$ in the consensus-ADMM iterates is a QP problem and can be solved by Theorem 4.

In the following example, we focus on the measurement ellipsoids obtained by Algorithm \ref{al4} using different numbers of the samples.

\textbf{Example 1.} Let $h(\hat v)=[\sqrt{\hat v_1^2+\hat v_2^2},\arctan({\hat v_2}/{\hat v_1})]^T$ for $\hat v=[\hat v_1,\hat v_2]^T$. Consider the problem \eqref{ADMM1} with $g(\hat x)=[1,-1]\hat x$ and ${\mathcal T}=\{r=h^{-1}(v_0-v), v\in{\mathcal E}{(0,\hat R)},v_0=h([100,100]^T)\}$, where $\hat R=diag(20^2,0.1^2)$. By generating $s$ random vectors from the ellipsoid ${\mathcal E{(0,\hat R)}}$, problem \eqref{ADMM1} is expressed as:
\begin{align}
&\min_{\hat{x},P} \quad \log \det(P)\nonumber \\
&s.t. \quad (r_i-\hat{x})^TP^{-1}(r_i-\hat{x})\leq 1,\nonumber \\
&\qquad\qquad [1,-1]\hat x=0,~i\in I,\label{example1sip1}
\end{align}
where $r_i\in{\mathcal T}$, $I=\{1,2,...,s\}$, and $\hat x=[\hat x_1,\hat x_2]^T$.
\begin{figure}
\begin{center}
\includegraphics[height=2.8cm]{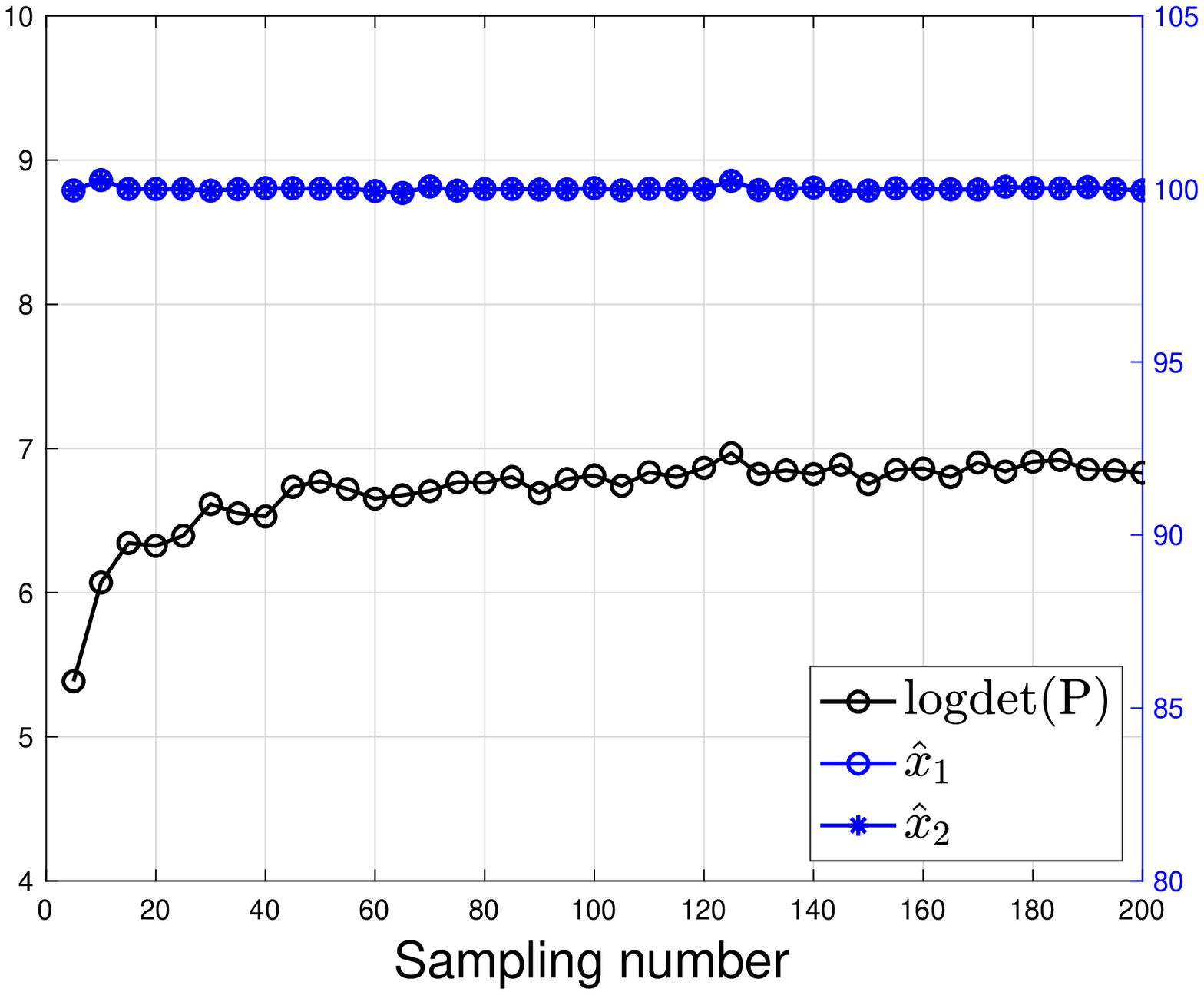} 
\includegraphics[height=2.8cm]{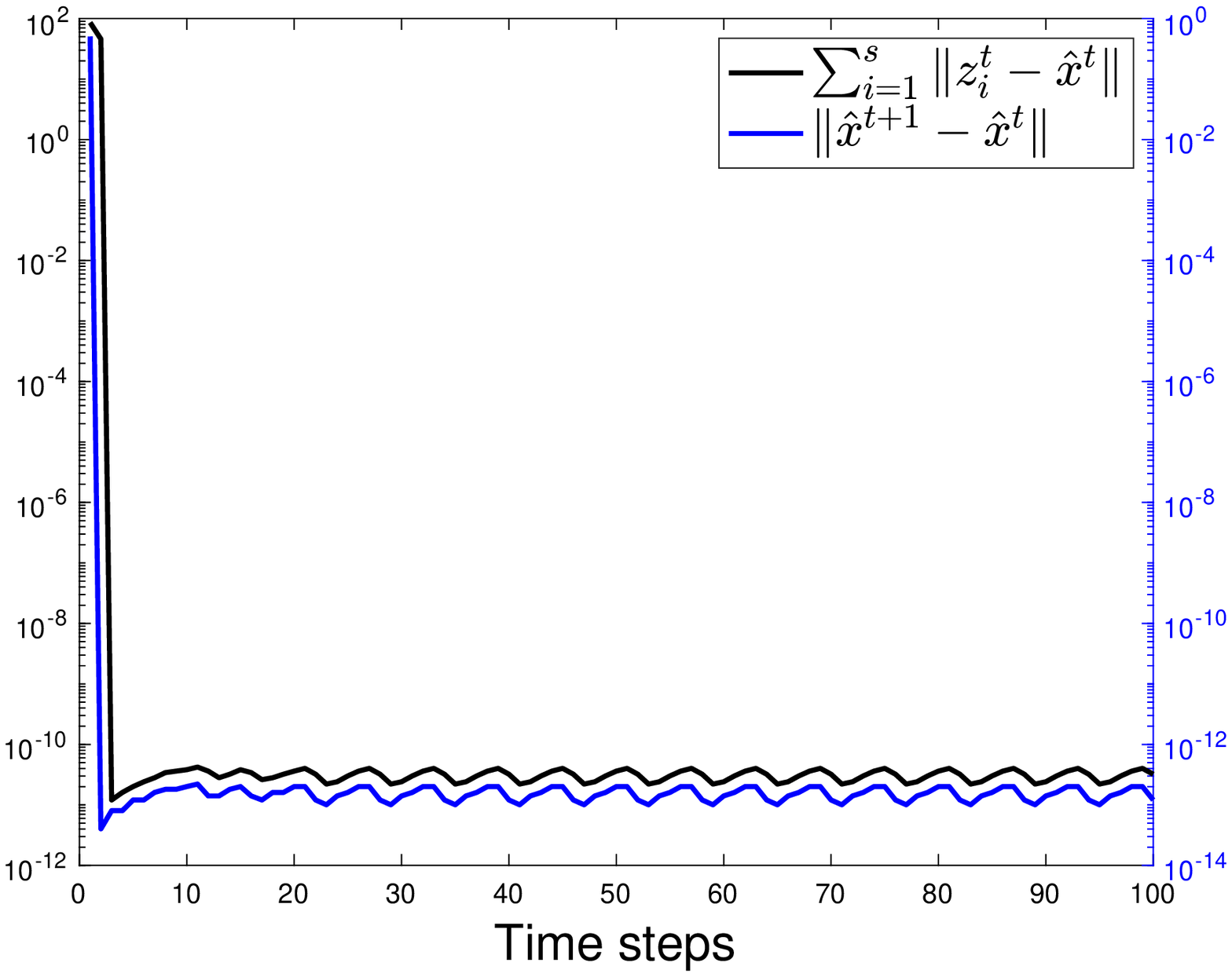} 
\caption{The logdet and the center of the estimation ellipsoids with different sampling number (left) and $\sum^s_{i=1}\|z_i^{t}-\hat{x}^{t}\|$ and $\|\hat{x}^{t+1}-\hat{x}^{t}\|$ versus iteration number with $s = 100$ (right).\label{figerror}}
\end{center}
\end{figure}
Fig. \ref{figerror} (left) plots the $\log \det(P)$ and the center of the bounding ellipsoid of set ${\mathcal T}$ obtained by solving problem \eqref{example1sip1} via Algorithm 3 for the different number of samples $s$. A total of 50 Monte Carlo runs are simulated. It shows that the logdet and the center can quickly converge to a stable value as the number of samples $s$ increases. Therefore, Algorithm \ref{al4} can get stable estimation ellipsoids with a small number of samples. Furthermore, $\sum^s_{i=1}\|z_i^{t}-\hat{x}^{t}\|$ and $\|\hat{x}^{t+1}-\hat{x}^{t}\|$ versus the number of iterations $t$ for one random problem instance with $s=100$ are plotted in Fig. \ref{figerror} (right). It shows that the differences between $z_i^{t}$ and $\hat{x}^{t}$, $\hat{x}^{t}$ and $\hat{x}^{t+1}$ rapidly decrease as the number of the algorithm iterations increases, i.e., \eqref{prop1}-\eqref{prop2} are satisfied indeed. In addition, the partial stopping criterion (the successive difference of $\hat x$ is smaller than a given tolerance) is also satisfied as the number of the iterations increases.

\subsection{Consensus-ADMM method versus semidefinite programming (SDP) method}
Consider the optimization problem \eqref{ADMM2} with a linear constraint ($g$ is a linear function), which can be rewritten by Schur complement \cite{boyd2004convex}:
\begin{align}
&\min_{\hat{x},P} \quad \log \det(P)\nonumber \\
&s.t. \quad \begin{bmatrix}
P&(r_i-\hat{x})\\(r_i-\hat{x})^T&1
\end{bmatrix}\succeq0, g(\hat{x})=0,~i\in I,\label{SDP1}
\end{align}
where $g(\hat x)=C\hat x-c$. The solution of the linear estimation constraint is $\hat{x}=C^\dagger c+U\gamma$, where $U=I-C^\dagger C$ is an orthogonal projection matrix in $\mathcal N(C)$ and $\gamma\in\mathbb{R}^n$ is arbitrary. Then the problem \eqref{SDP1} can be rewritten as an SDP problem by defining $H=P^{-1}$ and $\delta=P^{-1}U\gamma$:
\begin{align}
&\min_{\delta,H} \quad -\log \det(H)\nonumber \\
&s.t. \quad \begin{bmatrix}
H&H(r_i-C^\dagger c)-\delta\\(H(r_i-C^\dagger c)-\delta)^T&1
\end{bmatrix}\succeq0,\nonumber \\
&\qquad\quad i\in I.\label{SDP2}
\end{align}
In the problem \eqref{SDP2}, the dimension of decision variables is $M=\frac{n(n+1)}{2}+n$ and the dimension of the constraint matrix is $N=s(n+1)$. By using a primal-dual interior-point method to solve \eqref{SDP2}, the worst-case estimate requires $O(M^{2.75}s^{1.5})$ arithmetic operations to solve the problem with a given accuracy \cite{vandenberghe1996semidefinite,boyd1994linear}. Nevertheless, the consensus-ADMM method requires $O(n^3+sn^2)$ arithmetic operations to solve the problem \eqref{ADMM2} with a linear constraint. Specific solution steps and corresponding computational complexity are as
follows:
\begin{itemize}
  \item Update of $P$: Algorithm 2 is used to obtain the updated $\mu^*$. Each iteration of Algorithm 2 requires $O(n^2+(n+1)s)$ arithmetic operations, and calculating $P^*$ by equation \eqref{pb} requires $O(sn^2)$ arithmetic operations.
  \item Update of $z_i$: The arithmetic operations is $O(n)$ for each $i$ at each iteration by caching the eigen-decomposition of $P$ and using Newton’s method to solve the equation $g(\phi_i)=1$ \cite{huang2016consensus}. The eigen-decomposition of $P$ requires $O(n^3)$ arithmetic operations.
  \item Update of $\hat x$: Calculating $\hat{x}^*$ by equation \eqref{xxx} requires $O(n^3)$ arithmetic operations.
\end{itemize}
Therefore, consensus-ADMM method requires $O(n^3+sn^2)$ arithmetic operations to solve the problem \eqref{ADMM2} with a linear constraint, which is much lower than $O(M^{2.75}s^{1.5})$ arithmetic operations of the SDP method. 

For the problem \eqref{ADMM2} with a quadratic constraint, by similar analysis and using the Newton’s method to calculate the updated $\hat x^*$, the consensus-ADMM method also requires $O(n^3+sn^2)$ arithmetic operations to solve it. For the problem \eqref{ADMM2} with a general form constraint, the consensus-ADMM method requires $O(n^3+sn^2)$ arithmetic operations together with the computation of the constraint function and its gradient to solve it since the computational complexity of updating $\hat x$ depends on the specific form of the constraint function, and it costs computation of the constraint function and its gradient together with $O(n^3)$ arithmetic operations to obtain an updated $\hat{x}^*$ via a basic implementation of Newton’s method at each iteration \cite{bertsekas1997nonlinear}. Nevertheless, the SDP method is not suitable for solving the problem \eqref{ADMM2} with the nonlinear constraint.

In the following example, we compare the computing time for solving problem \eqref{SDP1} by using SDPT3 \cite{cvx} and the consensus-ADMM algorithm, respectively.

\textbf{Example 2.} Suppose there are $s$ samples $r_i\in\mathbb{R}^{\hat n}$ that are generated by a random generate ellipsoid $\mathcal E(0,\hat P)$. $C=[\hat C, -1, 0,...,0]\in\mathbb{R}^{1\times \hat n}$, $c=0$, and $\hat C$ is generated by the standard uniform distribution. The computing time of different methods for solving problem \eqref{SDP1} is given in Table \ref{tab1}. SDP indicates that the SDP problem \eqref{SDP2} is solved by SDPT3, C(S)-ADMM indicates sequential computation at the update of $z_i$ when solving the problem \eqref{SDP1} with the consensus-ADMM algorithm, and C(P)-ADMM indicates parallel computation at the update of $z_i$ when solving the problem \eqref{SDP1} with the consensus-ADMM algorithm. The results show that the consensus-ADMM algorithm takes less computing time than the SDP method. It is consistent with their computational complexity as $O(n^3+sn^2)$ and $O((\frac{n(n+1)}{2}+n)^{2.75}s^{1.5})$, respectively. In this example, the computation times of the consensus-ADMM algorithms are on the order of 0.1$\sim$1 sec. The alternative approach to further reduce the computation time of the proposed consensus-ADMM algorithm is using the hardware acceleration methods or selecting specific subproblem solvers according to the structure of the optimization problems.
\begin{table}[htbp]
    \centering 
    \caption{The mean computing time of SDPT3 and consensus-ADMM algorithms over 50 Monte Carlo runs.\label{tab1}}
    \renewcommand\arraystretch{1.2}
    \begin{tabular}{ccccc}\hline
      $\hat n$&$s$&SDP&C-ADMM(S)&C-ADMM(P)\\\hline
      5&50&1.0308&0.0909&0.0702\\
      5&100&1.6117&0.1382&0.1051\\
      5&200&2.7122&0.2102&0.1504\\
      5&400&5.2300&0.2823&0.1864\\
      5&800&10.6809&0.3623&0.1918\\
      5&1600&22.8754&0.5629&0.2400\\
      10&50&1.3213&0.1612&0.1027\\
      10&100&2.3113&0.3436&0.2563\\
      10&200&4.5404&0.4967&0.3778\\
      10&400&9.4156&0.6119&0.4334\\
     10&800&21.0134&0.7180&0.4075\\
      10&1600&53.0252&1.0270&0.4809\\\hline
    \end{tabular}
  \end{table}
\section{Simulation}\label{v}
In this section, we compare the performance of the consensus-ADMM-based method in Algorithm \ref{al3} with using the estimation constraints to both the prediction step and update step, the consensus-ADMM-based method in Algorithm \ref{al3} with only using the estimation constraints to the update step, dual set membership filter \cite{9435052} and the set membership filter with state constraints \cite{yang2009set1,yang2009set} by two numerical examples. To distinguish these methods, we abbreviate them as C-ADMM-SMF(EC), C-ADMM-SMF(MC), DSMF, and SMF-SC, respectively. Since SMF-SC considers the linear dynamic systems, we use the method in \cite{wang2019ellipsoidal} to linearize the nonlinear measurement function by performing Taylor expansion at the state prediction in \cite{9435052}.

\subsection{Linear estimation constraint}
In this example, we consider the problem of tracking a target in two dimensions, and the dynamic system is given by:
\begin{align}
x_{k+1}=\begin{bmatrix}
1&0&T&0\\
0&1&0&T\\ 
0&0&1&0\\
0&0&0&1
\end{bmatrix}x_k+w_k&,\label{51}\\
y_{k}=\begin{bmatrix}
\sqrt{(s^x_k-a)^2+(s^y_k-b)^2}\\
\arctan(\frac{s^y_k-b}{s^x_k-a})
\end{bmatrix}+&v_k,\label{52}
\end{align}
where $x_k:=(s^x_k, s^y_k, \sigma^x_k, \sigma^y_k)$ and $T=1$ is the sampling time. The noises are assumed to be restricted in $w_k\in {\mathcal E}(0,Q_k)$ and $v_k\in {\mathcal E}(0,R_k)$, respectively, where $Q_k=10I_{4\times4}$ and $R_k=diag(20^2,0.1^2)$. The estimation constraint on the above system is $C\hat x_k=c$, where $C=\begin{bmatrix}
2&-1&0&0\\
0&0&2&-1
\end{bmatrix}$ and $c=[0,0]^T$.

In the simulation, the initial state $x_0=[0, 0, 25, 50]^T$ and $a=15000$, $b=0$. Assume the initial estimate $\hat{x}_0$ of the state is a random disturbance around the initial state, $P_0 = 100^2I_{4\times4}$, and $w_k = (\sin(k\pi/2), 2\sin(k\pi/2), \sin(k\pi/4), 2\sin(k\pi/4))$. A total of 100 Monte Carlo runs are simulated. The root mean square error (RMSE) is defined as $RMSE_k=\sqrt{\frac{1}{L}\sum^L_{i=1}(\hat x^i_k-x^i_k)^2}$, where $\hat x^i_k$ and $x^i_k$ are the state estimation and true state at the $k$-th time step and $i$-th Monte Carlo, respectively. $L$ is the number of the Monte Carlo runs.

\begin{figure}[htbp]
\begin{center}
\includegraphics[height=4cm]{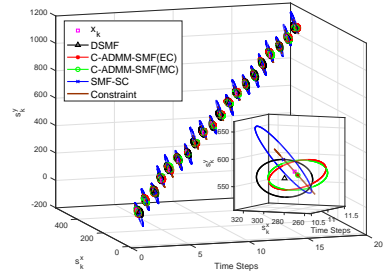}  
\caption{The true state $x_k$ and the state estimates.}
\label{fig001}                                 
\end{center}                                
\end{figure}
\begin{figure}[htbp]
    \centering
    \includegraphics[height=2.8cm]{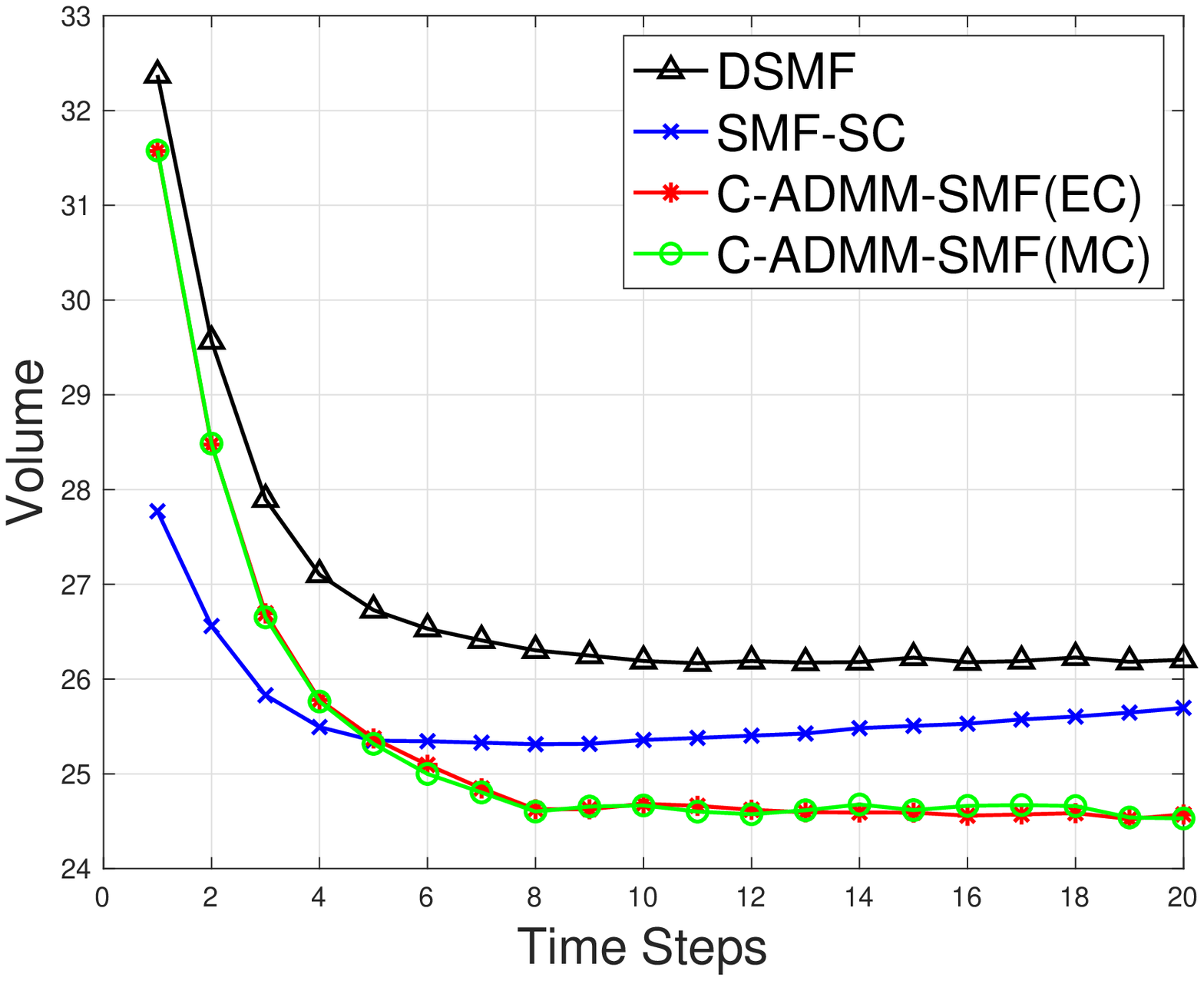}
    \includegraphics[height=2.8cm]{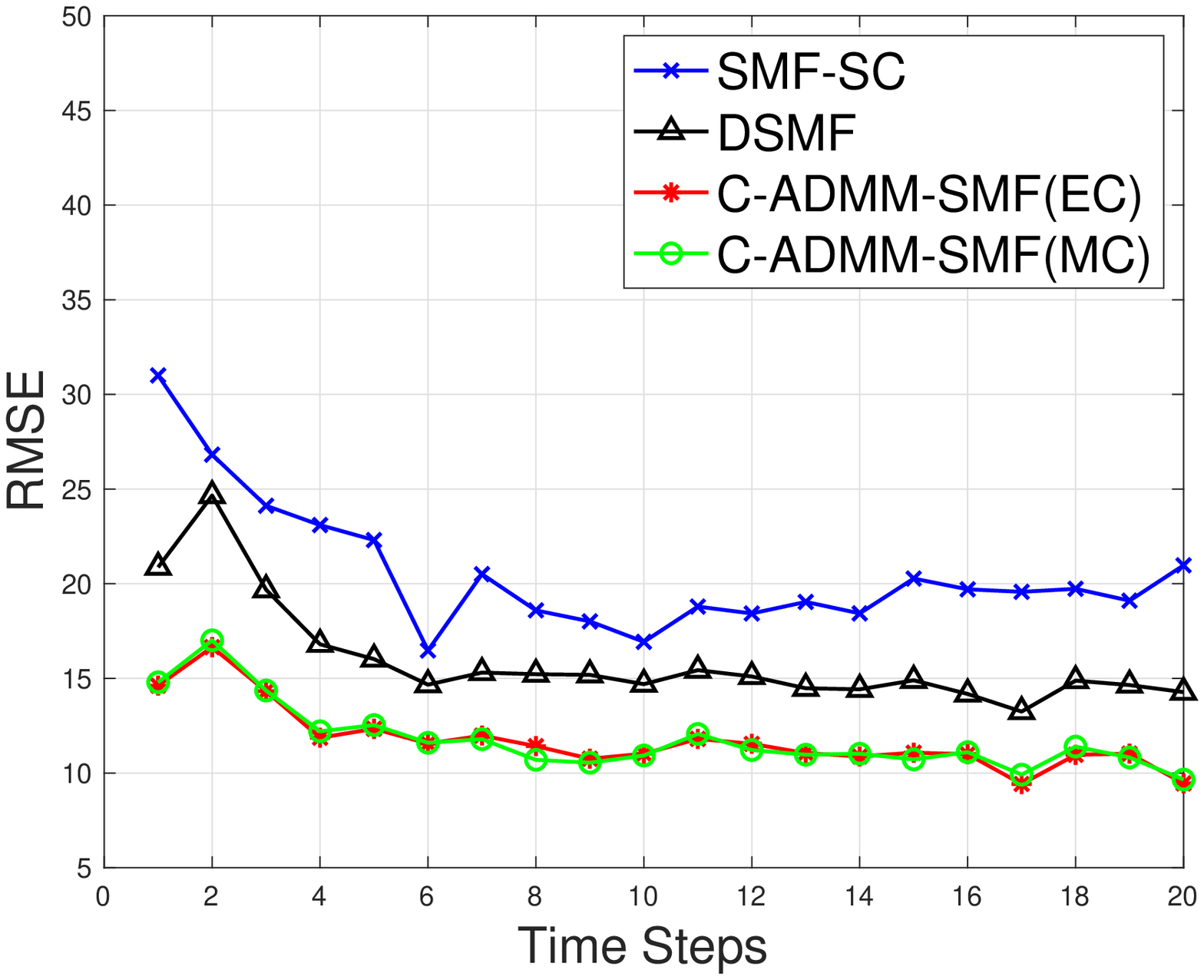}  
\caption{The volume (logdet) of the state bounding ellipsoids (left) and the RMSE of the state estimates (right).}
\label{fig3} 
\end{figure}
\begin{figure}[htbp]
    \centering
   \includegraphics[height=2.8cm]{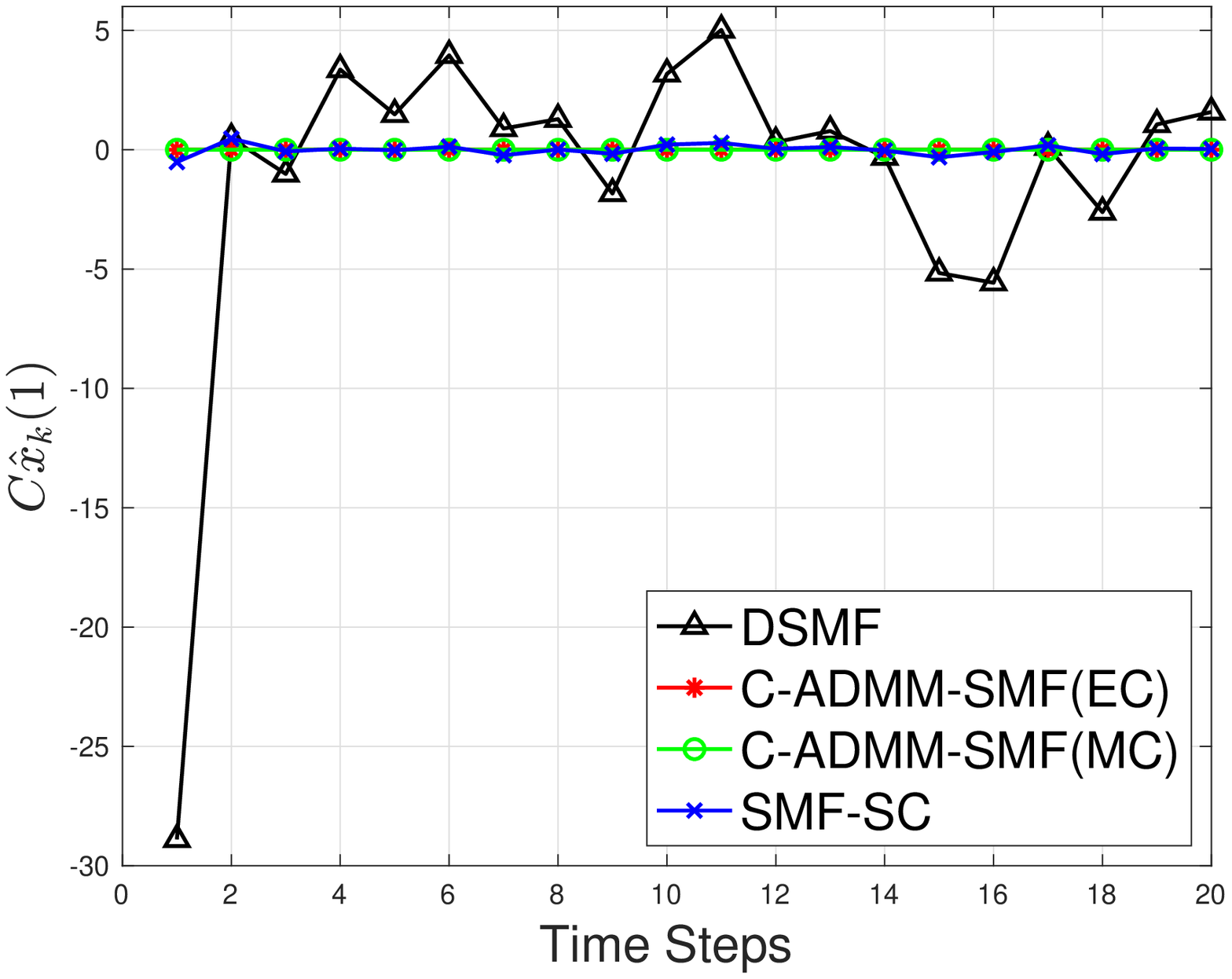}  
\includegraphics[height=2.8cm]{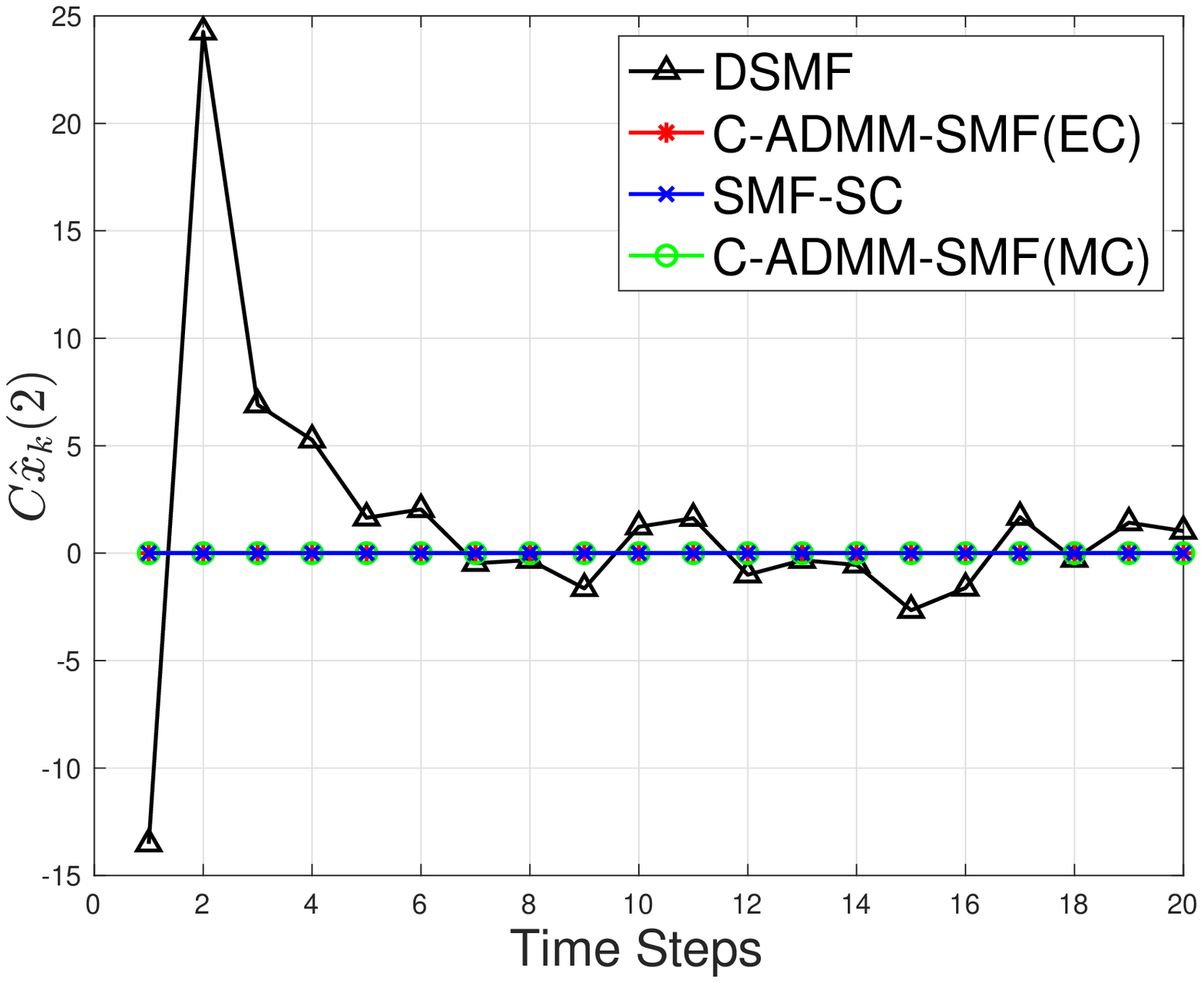}  
\caption{Satisfaction of the constraints.}
\label{fig5} 
\end{figure}

Fig. \ref{fig001} shows the trajectories of the true state, the state estimates, and the state bounding ellipsoids obtained by different methods. Fig. \ref{fig3} plots the logdet of the shape matrix of the state bounding ellipsoids and the state estimation error versus time steps obtained by different methods. Figs. \ref{fig001}-\ref{fig3} show that the size of the state bounding ellipsoids obtained by C-ADMM-SMF(EC) and C-ADMM-SMF(MC) is smaller than that of DSMF and SMF-SC and the RMSE of C-ADMM-SMF(EC) and C-ADMM-SMF(MC) is less than that of DSMF and SMF-SC. Consistent with the results in Fig. \ref{fig001}, Fig. \ref{fig5} shows that the state estimation provided by DSMF and SMF-SC do not always lie on the constraint, while C-ADMM-SMF(EC) and C-ADMM-SMF(MC) do produce estimates that satisfy the constraint. The reason is that DSMF does not use the constraint information and C-ADMM-SMF(EC) and C-ADMM-SMF(MC) use the SIP approach to transform the nonlinear system into a linear one to obtain a more accurate estimation ellipsoid and use the estimation constraint information to produce the estimation that satisfies the constraint. In addition, the performance of C-ADMM-SMF(EC) and C-ADMM-SMF(MC) is similar. The reason may be that C-ADMM-SMF(EC) produces a more accurate predicted ellipsoid, so the state bounding ellipsoid determined by C-ADMM-SMF(EC) is similar to that of C-ADMM-SMF(MC).

\subsection{Nonlinear estimation constraint}
In the case of nonlinear estimation constraints, we consider the dynamic system of the following form:
\begin{align}
x_{k+1}=\begin{bmatrix}
1&0&\frac{\sin wT}{w}&-\frac{1-\cos wT}{w}\\
0&1&\frac{1-\cos wT}{w}&\frac{\sin wT}{w}\\ 
0&0&\cos wT&-\sin wT\\
0&0&\sin wT&\cos wT
\end{bmatrix}x_k+w_k,\label{}
\end{align}
where $x_k:=(s^x_k, s^y_k, \sigma^x_k, \sigma^y_k)$, the measurement function and the settings for noises $w_k$ and $v_k$ are the same as in the linear case. The state estimation constraints on this dynamic system are $\hat x_k^TC_1\hat x_k=c_1^2$ and $\hat x_k^TC_2\hat x_k=c_2^2$, where $C_1=\begin{bmatrix}
I_{2\times 2}&0_{2\times 2}\\0_{2\times 2}&0_{2\times 2}
\end{bmatrix}$, $C_2=\begin{bmatrix}
0_{2\times 2}&0_{2\times 2}\\0_{2\times 2}&I_{2\times 2}
\end{bmatrix}$, $c_1=10$, and $c_2=0.05$.

In this simulation, the sampling time is $T=1$,
$a=-20000$, $b=5000$, and the initial state is $x_0=(0, 10, -0.05, 0)$. The settings of the initial state bounding ellipsoid and the total Monte Carlo runs are the same as in the linear case.

\begin{figure}
\begin{center}
\includegraphics[height=4cm]{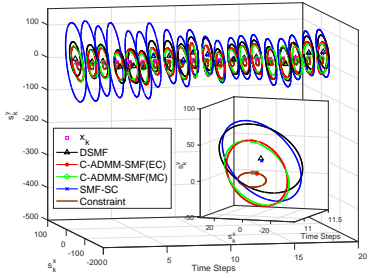}  
\caption{The true state $x_k$ and the state estimates.}
\label{fig002}                                 
\end{center}                                
\end{figure}
\begin{figure}[htbp]
  \centering
    \includegraphics[height=2.8cm]{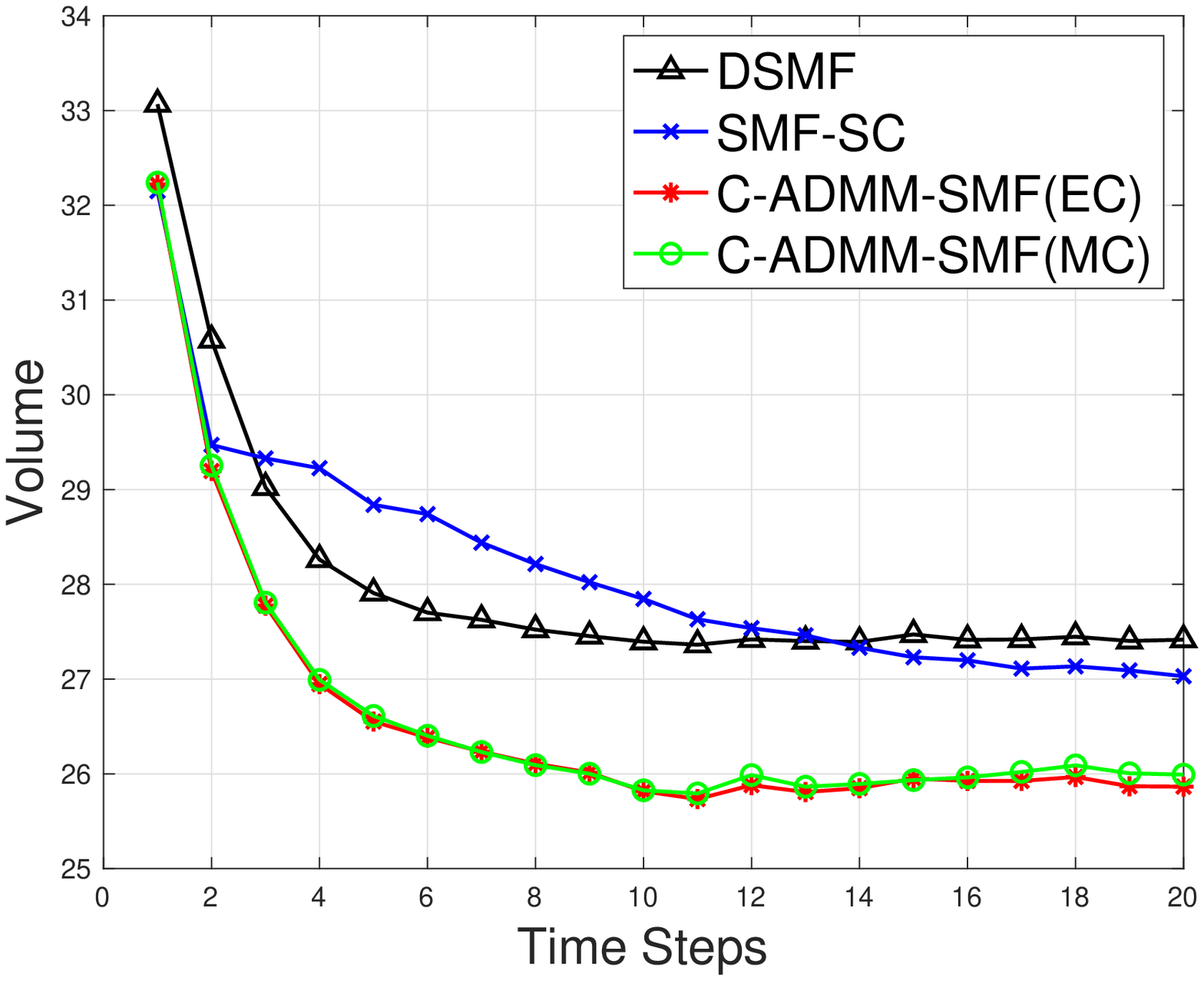}
    \includegraphics[height=2.8cm]{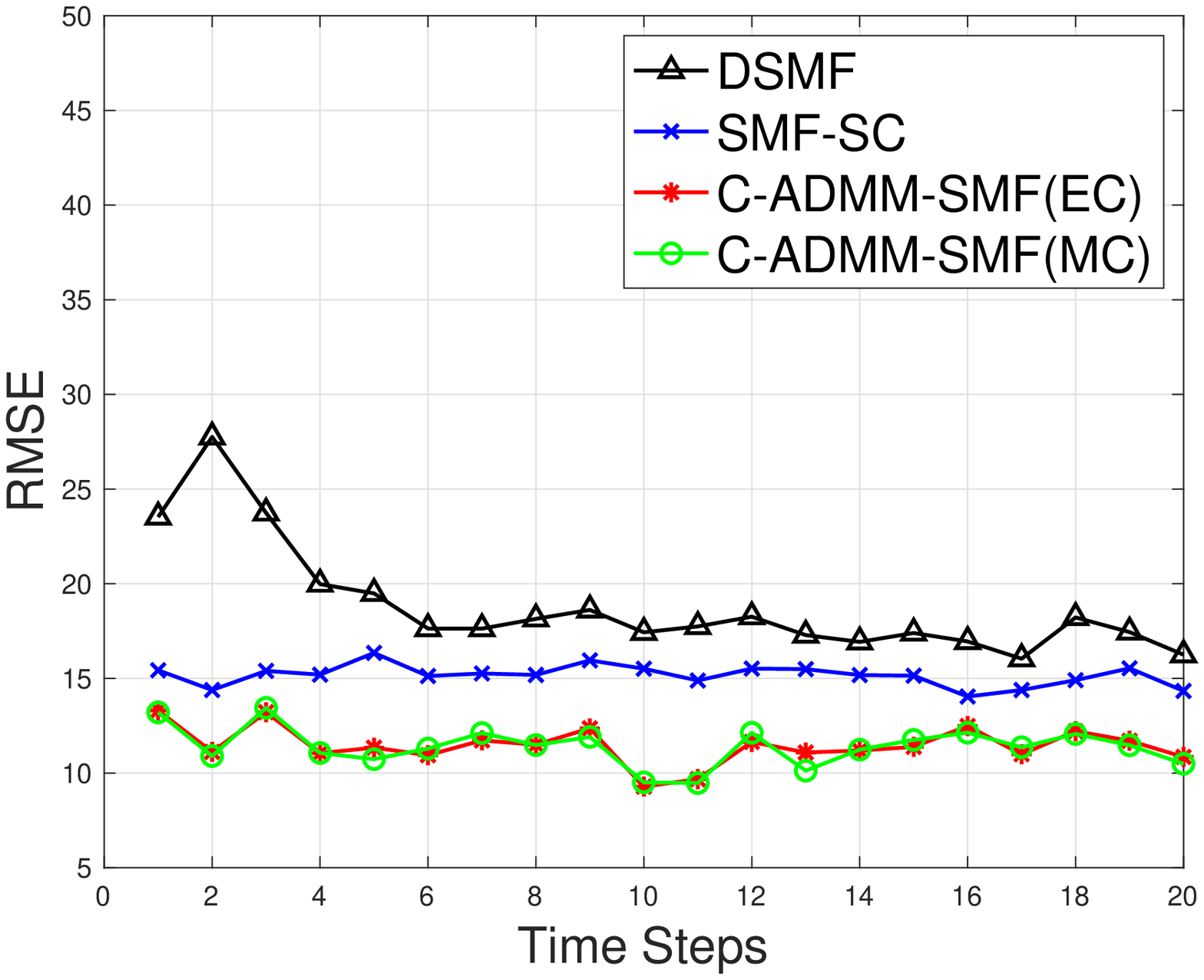}  
\caption{The volume (logdet) of the state bounding ellipsoids (left) and the RMSE of the state estimates (right).}
\label{fign3} 
\end{figure}
\begin{figure}[htbp]
    \centering
   \includegraphics[height=2.8cm]{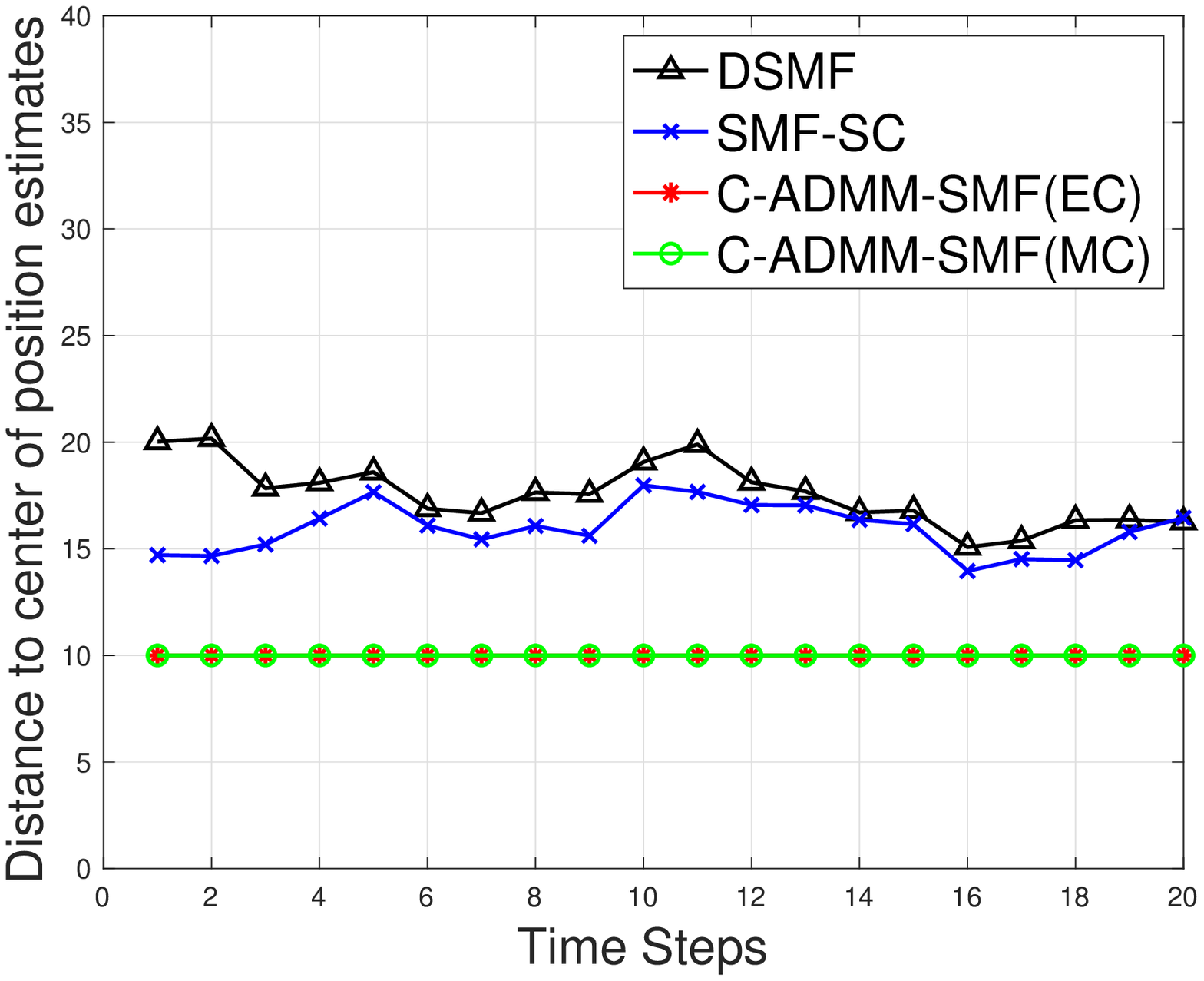}  
\includegraphics[height=2.8cm]{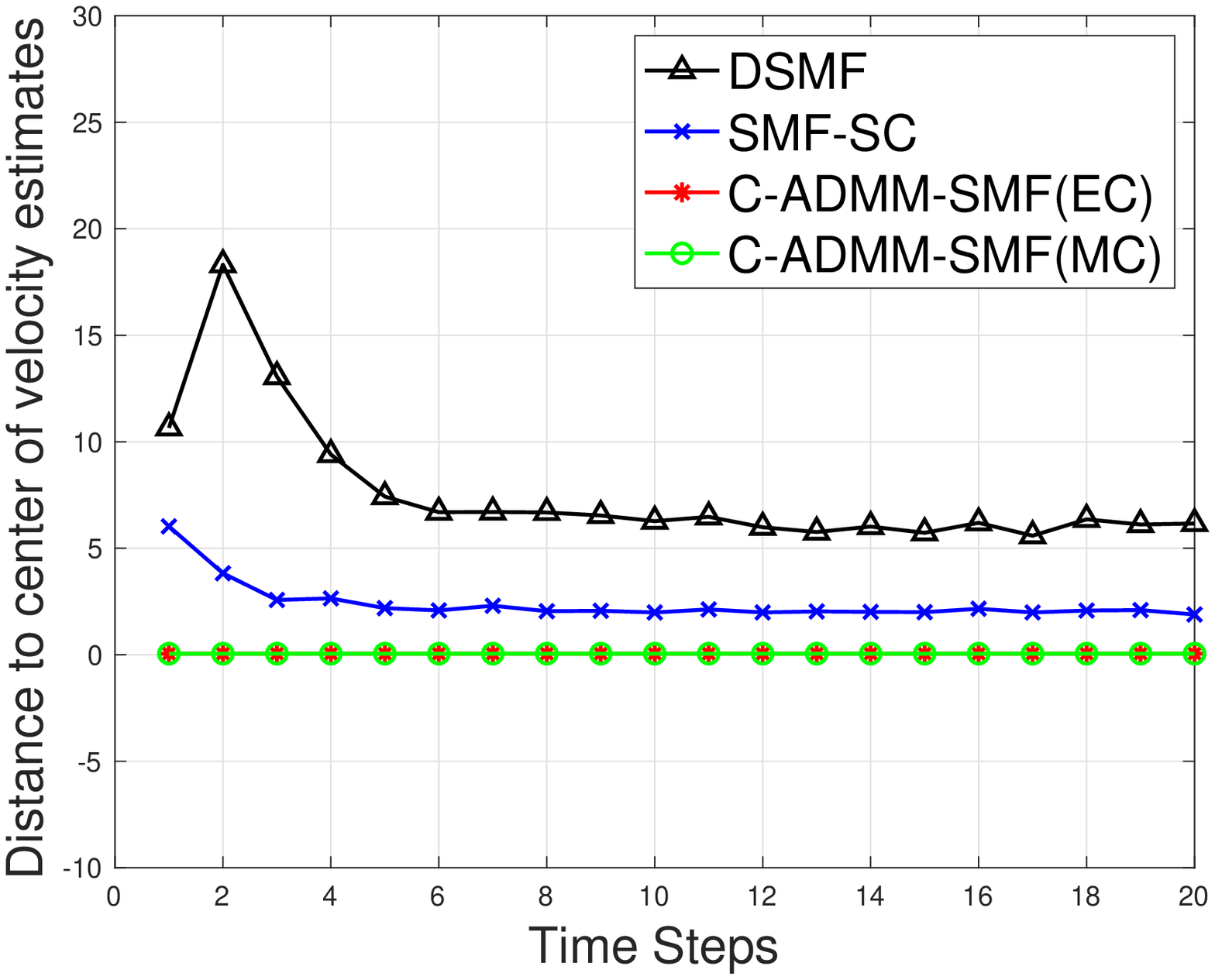}  
\caption{Satisfaction of the constraints.}
\label{fign5} 
\end{figure}

The simulation results are similar to the case of the linear estimation constraint. Figs. \ref{fig002}-\ref{fign3} show the trajectories of the true state, the state estimates, and the state bounding ellipsoids, logdet of the shape matrix of the estimation ellipsoids, and the RMSE of the state estimation versus time steps by C-ADMM-SMF(EC), C-ADMM-SMF(MC), DSMF, and SMF-SC, respectively. They show that C-ADMM-SMF(EC) and C-ADMM-SMF(MC) perform better than DSMF and SMF-SC, and the performance of C-ADMM-SMF(EC) and C-ADMM-SMF(MC) are similar. Fig. \ref{fign5} presents the distance between the state estimations and the center point of the constraints. Consistent with the results in Fig. \ref{fig002}, which also shows that the state estimation provided by DSMF and SMF-SC do not always satisfy the constraint, while C-ADMM-SMF(EC) and C-ADMM-SMF(MC) do. Similar to the linear constraint case, the reason is also that C-ADMM-SMF(EC) and C-ADMM-SMF(MC) use the SIP approach to transform the nonlinear system into a linear one instead of linearizing the nonlinear system to obtain a more accurate estimation ellipsoid.

\section{Conclusion}\label{vi}
This paper considered the problem of state estimation for nonlinear dynamic systems with unknown but bounded noises and state estimation constraints. We developed a recursive set membership algorithm to compute the state bounding ellipsoid that includes the prediction and measurement update steps, with the center of the ellipsoid satisfies the constraint. The nonlinear dynamic system is transformed into a linear system by solving the SIP problems instead of linearizing the nonlinear functions, which allows us to obtain a tighter state bounding ellipsoid. In addition, a consensus-ADMM-based algorithm is proposed to solve the SIP problems and each iteration of the algorithm can be solved efficiently. Finally, typical numerical examples have demonstrated the effectiveness of the proposed method. Future research directions may include the state estimation problem for more general dynamic systems (e.g., maintaining the state estimation performance when the measurement function is not invertible or the dynamic systems are nonlinear in the noises), the investigation of the more efficient algorithms for solving the nonconvex SIP problems (e.g., exploring algorithms with lower computational complexity), and the generalization of the proposed filter to multi-sensor fusion.

\bibliographystyle{unsrt}        
\bibliography{ref}           

\appendix

\section{The proof of Lemma 1} \label{app0}
\textbf{Proof:} According to Theorem 4.1 in \cite{durieu2001multi} and (\ref{pre1}), the center of the predicted ellipsoid ${\mathcal E}_{k+1|k}$ is given by $\hat{x}_{k+1|k}=\hat{x}_{f_k}+0$, thus we obtain the result. \hfill{$\square $}

\section{The proof of Theorem 1} \label{appa}  
\textbf{Proof:} Define $H=P^{-1}$, then (\ref{ADMM2}) can be rewritten as:
\begin{align}
&\min_{\hat x, H} \quad -\log \det(H)\nonumber \\
&s.t. \quad (r_i-\hat{x})^TH(r_i-\hat{x})\leq 1,\nonumber \\
&\qquad\qquad g(\hat{x})=0,\nonumber \\
&\qquad\qquad i\in I.\label{ADMM3}
\end{align}
A KKT point $\{H^*,\hat x^*\}$ of (\ref{ADMM3}) and the corresponding dual variables $\varsigma_i^*$ and $\varkappa^*$ satisfy that:
\begin{align}
-(H^*)^{-1}+\sum_{i\in I}\varsigma_i^*(r_i-\hat x^*)(r_i-\hat x^*)^T=0,\nonumber\\
\varkappa^*\nabla g(\hat x^*)-2\sum_{i\in I}\varsigma_i^*H^*(r-\hat x^*)=0,\nonumber\\
\varsigma_i^*\geq0,\nonumber\\
(r_i-\hat{x}^*)^TH^*(r_i-\hat{x}^*)\leq 1,\nonumber\\
 g(\hat{x}^*)=0,\nonumber\\
 \varsigma_i^*((r_i-\hat{x}^*)^TH^*(r_i-\hat{x}^*)-1)=0,\nonumber\\
\forall i\in I.\nonumber
\end{align}
Let the superscript $t$ denote the point obtained at iteration $t$. At iteration $t+1$, let $H=P^{-1}$ in the optimization (\ref{pz1}) and assume that each $z_i^{t}$ is well defined, from (\ref{pz1}) - (\ref{ADMM11111}) we have:
\begin{align}
-(H^{t+1})^{-1}+\sum_{i\in I}\varpi_i^{t+1}(r_i-z_i^{t+1})(r_i-z_i^{t+1})^T=0,\\
-\lambda_i^{t}-\rho(\hat x^{t}-z_i^{t+1})-2\varpi_i^{t+1}H^{t+1}(r_i-z_i^{t+1})=0,\label{zi}\\
(r_i-z_i^{t+1})^TH^{t+1}(r_i-z_i^{t+1})\leq 1,\\
\varpi_i^{t+1}\geq 0\\
\varpi_i^{t+1}((r_i-z_i^{t+1})^TH^{t+1}(r_i-z_i^{t+1})-1)=0\\
\sum_{i\in I}\lambda_i^{t}+\rho\sum_{i\in I}(\hat x^{t+1}-z_i^{t+1})+\psi^{t+1}\nabla g(\hat x^{t+1})=0,\label{xx1}\\
g(\hat{x}^{t+1})=0,\\
\forall i\in I,\nonumber
\end{align}
where $\varpi_i^{t+1}$ and $\psi^{t+1}$ are dual variables of (\ref{pz1}) and (\ref{x0}), respectively.
By adding equation (\ref{zi}) from $i=1$ to $i=s$, and use the fact that $\lambda_i^{t+1}=\lambda_i^{t}+\rho(\hat{x}^{t+1}-z_i^{t+1})$ we have:
\begin{align}
-\sum_{i\in I}\lambda_i^{t+1}&+\rho s(\hat x^{t+1}-\hat x^{t})\nonumber\\
&-2\sum_{i\in I}\varpi_i^{t+1}H^{t+1}(r_i-z_i^{t+1})=0.
\end{align}
Furthermore, based on the assumption and (\ref{xx1}), we can get:
\begin{align}
\psi^{t+1}\nabla g(\hat x^{t+1})-2\sum_{i\in I}\varpi_i^{t+1}H^{t+1}(r_i-z_i^{t+1})=0.
\end{align}
Let $\varpi_i^{t+1}=\varsigma_i^*$ and $\psi^{t+1}=\varkappa^*$, then the rest of the KKT conditions can be guaranteed by the assumptions.\hfill{$\square $}

\section{The proof of Theorem 2} \label{appb} 
\textbf{Proof:} Define $H=nP^{-1}$ and $\alpha_i=r_i-z_i$, (\ref{P111}) can be rewritten as:
\begin{align}
 \min_{H}\quad &-\log \det(H),\nonumber\\
s.t. \quad &\alpha_i^TH\alpha_i\leq n,\nonumber\\
&\quad i\in I.\label{p111}
\end{align}
The Lagrangian function of (\ref{p111}) is:
\begin{align}
&L(H,\mu)=-\log \det(H)+\sum_{i\in I}\mu_i(\alpha_i^TH\alpha_i-n),\nonumber
\end{align}
and the dual problem of (\ref{p111}) is
\begin{align}
\max_{\mu}\min_{H}&\quad L(H,\mu)\nonumber\\
s.t.&\quad \mu\geq 0.\label{dual}
\end{align} 
The inner minimum in (\ref{dual}) is achieved at $H^*$ if and only if:
\begin{align}
0=\nabla_HL(H^*,\mu)=-{H^*}^{-1}+\sum_{i\in I}\mu_i\alpha_i\alpha_i^T.\label{H}
\end{align}
From (\ref{H}) we have:
\begin{align}
L(H^*,\mu)=\log \det&(\sum_{i\in I}\mu_i\alpha_i\alpha_i^T)-n\sum_{i\in I}\mu_i\nonumber\\
&+\sum_{i\in I}\mu_i\alpha_i^T(\sum_{i\in I}\mu_i\alpha_i\alpha_i^T)^{-1}\alpha_i,
\end{align}
where
\begin{align}
&\sum_{i\in I}\mu_i\alpha_i^T(\sum_{i\in I}\mu_i\alpha_i\alpha_i^T)^{-1}\alpha_i\nonumber\\
&=\sum_{i\in I}tr(\mu_i\alpha_i^T(\sum_{i\in I}\mu_i\alpha_i^T)^{-1}\alpha_i)\nonumber\\
&=\sum_{i\in I}tr((\sum_{i\in I}\mu_i\alpha_i\alpha_i^T)^{-1}\mu_i\alpha_i\alpha_i^T)\nonumber\\
&=tr((\sum_{i\in I}\mu_i\alpha_i\alpha_i^T)^{-1}(\sum_{i\in I}\mu_i\alpha_i\alpha_i^T))\nonumber\\
&=n.\nonumber
\end{align}
Thus, the dual problem \eqref{dual} is:
\begin{align}
\max_\mu&\quad \log \det(\sum_{i\in I}\mu_i\alpha_i\alpha_i^T)-n\sum_{i\in I}\mu_i+n\nonumber\\
s.t.& \quad\mu\geq 0.\label{dual111}
\end{align} 

For any $\hat{\mu}$ can be written as $t\mu$, where $t$ is nonnegative and $\sum_{i\in I}\mu_i=1, \mu\geq 0$, then
\begin{align}
 &\log \det(\sum_{i\in I}\hat\mu_i\alpha_i\alpha_i^T)-n\sum_{i\in I}\hat\mu_i+n\nonumber\\&
=\log \det(t\cdot\sum_{i\in I}\mu_i\alpha_i\alpha_i^T)-nt+n\nonumber\\&
=\log \det(\sum_{i\in I}\mu_i\alpha_i\alpha_i^T)+n\cdot logt-nt+n,
\end{align}
and this is maximized by choosing $t=1$. Thus, (\ref{dual111}) becomes:
\begin{align}
\max_\mu&\quad \log \det(\sum_{i\in I}\mu_i\alpha_i\alpha_i^T)\nonumber\\
s.t.&\quad\sum_{i\in I}\mu_i=1, \mu\geq 0.\label{p2}
\end{align} 
Thus, the optimal solution of (\ref{p111}) is
\begin{align}
{H}^*=(\sum_{i\in I}\mu_i^*\alpha_i\alpha_i^T)^{-1},
\end{align} 
where $\mu^*$ is the optimal solution of (\ref{p2}), which means 
\begin{align}
{P}^*=n\sum_{i\in I}\mu_i^*(r_i-z_i)(r_i-z_i)^T
\end{align}
is the optimal solution of (\ref{P111}). \hfill{$\square $}

\section{The proof of Theorem 3} \label{appc}   
\textbf{Proof:}
By defining $\chi_i=z_i-r_i$, (\ref{ADMM_pz}) can be written as
\begin{align}
\min_{\chi_i}& \|\chi_i-(\hat{x}-r_i+\frac{1}{\rho}\lambda_i)\|^2_2,\nonumber\\
&s.t. \quad \chi_i^TP^{-1}\chi_i\leq 1.\label{z1}
\end{align}

Clearly, the optimization problem without the constraint $\chi_i^TP^{-1}\chi_i\leq 1$ can be minimized as $\chi_i=\hat{x}-r_i+\frac{1}{\rho}\lambda_i$. Therefore if:
\begin{align}
(\hat{x}-r_i+\frac{1}{\rho}\lambda_i)^TP^{-1}(\hat{x}-r_i+\frac{1}{\rho}\lambda_i)\leq 1,
\end{align}
we have $\chi_i^*=\hat{x}-r_i+\frac{1}{\rho}\lambda_i$ is the optimal solution. If not, the optimal solution $\chi_i^*$ must satisfy ${\chi_i^*}^TP^{-1}\chi_i=1$. Thus (\ref{z1}) can be reformulated as:
\begin{align}
\min_{\chi_i}& \|\chi_i-(\hat{x}-r_i+\frac{1}{\rho}\lambda_i)\|^2_2\nonumber\\
&s.t. \quad \chi_i^TP^{-1}\chi_i=1,\label{z2}
\end{align}
and the Lagrangian function of (\ref{z2}) is:
\begin{align}
L(\chi_i,\phi_i)=\|\chi_i-(\hat{x}&-r_i+\frac{1}{\rho}\lambda_i)\|^2_2\nonumber\\&+\phi_i(\chi_i^TP^{-1}\chi_i-1).
\end{align}
Denote $\chi_i^*$ as the optimal solution of (\ref{z2}), based on the KKT optimality conditions, we have:
\begin{align}
&\nabla_{\chi_i}L(\chi_i^*,\phi_i)=0,\label{z3}
\end{align}
for some $\phi_i\geq0$. From (\ref{z3}) we have:
\begin{align}
\chi_i=(I+&\phi_iP^{-1})^{-1}(\hat{x}-r_i+\frac{1}{\rho}\lambda_i).\label{z4}
\end{align}
Let $P=E\Gamma E^T$ be the eigenvalue decomposition of $P$ and plugging (\ref{z4}) back into the equality constraint $\chi_i^TP^{-1}\chi_i=1$, we have:
\begin{align}
\sum^n_{j=1}\frac{d_j^2}{(\phi_i+\beta_j)^2}=1,
\end{align}
where $\Gamma=diag\{\beta_1,\beta_2,...,\beta_n\}$, and $d=E^TP^{\frac{1}{2}}(\hat{x}-r_i+\frac{1}{\rho}\lambda_i)$. Let
\begin{align}
g(\phi_i)=\sum^n_{j=1}\frac{d_j^2}{(\phi_i+\beta_j)^2},
\end{align}
we have:
\begin{align}
g^{'}(\phi_i)=-2\sum^n_{j=1}\frac{d_i^2}{(\phi_i+\beta_i)^3}<0.
\end{align}
Because of $g(0)>1$ and $g(\phi_i)$ monotonically decreases to zero as $\phi_i\rightarrow\infty$. Therefore the equation $g(\phi_i)=1$ has exactly one nonnegative solution $\phi_i^*$, and the optimal solution is $\chi_i^*=(I+\phi_i^*P^{-1})^{-1}(\hat{x}-r_i+\frac{1}{\rho}\lambda_i)$. Let $z_i^*=\chi_i^*+r_i$, we can obtain the final result.\hfill{$\square $}

\section{The proof of Theorem 4} \label{appd}  
\textbf{Proof:}
Consider the linear equality constraint $C\hat x-c=0$ with the solution $\hat{x}=C^\dagger c+U\gamma$, where $\gamma\in\mathbb{R}^n$ is arbitrary. Then, the optimization problem (\ref{123}) of updating $\hat{x}$ is transformed into:
\begin{align}
\min_{\gamma}\sum_{i\in I}\lambda_i^T(C^\dagger c+U\gamma-z_i)+\frac{\rho}{2}\sum_{i\in I}\|C^\dagger c+U\gamma-z_i\|^2,\label{x2}
\end{align}
which can be further simplified as:
\begin{align}
\min_{\gamma}\|U\gamma+C^\dagger c+\frac{1}{s}\sum_{i\in I}(\frac{\lambda_{i}}{\rho}-z_{i})\|^2_2,\label{x3}
\end{align}
where the solution is:
\begin{align}
\gamma=U^\dagger (\frac{1}{s}\sum_{i\in I}(z_{i}-\frac{\lambda_{i}}{\rho})-C^\dagger c).\label{x4}
\end{align}
Substituting (\ref{x4}) into $\hat{x}=C^\dagger c+U\gamma$ yields the result.\hfill{$\square $}}

\end{document}